\documentclass[showpacs,floatfix,amsmath,amssymb,floatfix]{revtex4}
\usepackage{graphicx,color,dcolumn,booktabs,bm}
\usepackage{longtable,lscape}
\usepackage{txfonts}
\usepackage{overpic}
\usepackage{amssymb}
\usepackage{indentfirst}
\usepackage{CJK}\usepackage{graphics}\usepackage{subfigure}\usepackage{tabularx}\usepackage{booktabs}
\usepackage{amsmath}\usepackage{epsfig}\usepackage{indentfirst}\usepackage{flafter}\usepackage{multirow}\usepackage{threeparttable}
\usepackage{amssymb}\usepackage{longtable}
\usepackage{booktabs}
\usepackage{array}
\usepackage{feynmf}   
\usepackage{slashed}  
\usepackage{cases}
\usepackage{color}
\usepackage{multirow}
\usepackage{graphicx,color,dcolumn,booktabs,bm}

\graphicspath{{Figures/}} %

\graphicspath{{Figures/}} %
\usepackage[colorlinks, citecolor=black,anchorcolor=black,menucolor=black, linkcolor=black,filecolor=black,runcolor=black,urlcolor=blue,frenchlinks=black]{hyperref}
\usepackage{cases}
\usepackage{CJK}\usepackage{graphics}\usepackage{subfigure}\usepackage{tabularx}\usepackage{booktabs}
\usepackage{amsmath}\usepackage{epsfig}\usepackage{indentfirst}\usepackage{flafter}\usepackage{multirow}\usepackage{threeparttable}
\usepackage{amssymb}\usepackage{longtable}
\usepackage{booktabs}
\usepackage{array}

\begin{document}

\title{Equivalent of the elliptic function solutions
of nonlinear differential equations}
\author{Dong-Hua Luo$^{1}$}
\author{Cheng-Qun Pang$^{2,3}$}\email{pangchq13@lzu.edu.cn}
\affiliation{$^1$Shandong Communication Media College ,~Shandong Jinan 25020,~China\\
$^2$Department of Physics,Liupanshui Normal University ~~Liupanshui,~~~553004,~China\\
$^3$School of Physical Science and Technology,~Lanzhou University,~Lanzhou 730000,~China}

\begin{abstract}
In this paper, we point out that many Jacobi elliptic function solutions to non-linear differential equation(NDE)  can be transformed each other via the modulus and phase transformation of Jacobi elliptic function.~Therefore these solutions are equivalent. We  investigate  equivalent of the  Jacobi elliptic function solutions of Korteweg-de Vries (mKdV) equation as a example.
\end{abstract}

 \maketitle

\section{Introduction}\label{sec1}
In recent decades, nonlinear differential equation(NDE)  become research focus,~so the methods for solving them have been developed rapidly.~Among them,~a series  Jacobi elliptic function solutions of nonlinear equations  were given by the means of  the Jacobi elliptic function expansion method$^{[1]}$ and its modified form$^{[2,3]}$ .~In fact,~many forms of Jacobi elliptic function solutions  can be obtained from one of them using the transformation of the Jacobi elliptic function, therefore, they are equivalent.~Further more, Jacobi elliptic function expansion method and its related correction method to reduce the number of expanded form,~or expansion to certain restrictions,~such as in Ref.[4] they have mentioned a few  transformations of the Jacobi elliptic function,~and used them to simplify the Jacobi elliptic function expansion.~In Ref.[5],~they pointed out that Jacobi elliptic transformation  can be use to obtain other Jacobi elliptic functional solutions of NDE.~In this paper, we point out that we should test wether  the solutions are equivalent and we give a more systematic and comprehensive shown for transformation of the Jacobi elliptic function which can help us to test.~Basing on any one of the known NDE's Jacobi elliptic function solutions,~you can get other Jacobi elliptic function solutions by the transformation, of cause they are equivalent solutions.~Finally,~we investigate  Jacobi elliptic function solutions to mKdV equation as a example.
\par
The paper is organized as follows: In section II, we review the modulus and phase transformation of Jacobi elliptic function and give some new forms of them. A treatment for the Jacobi elliptic function solutions to mKdV equation as an example is shown in section III. In the last section, we give a brief conclusion.

\section{Jacobi elliptic transformation}\label{sec2}

The definition and some properties of the Jacobi elliptic function can see appendix.
Assuming $\text{sn}(x+\phi | \lambda )$ can become S($x)$ by the module and phase transformation,~at the same time $\text{cn}(x+\phi | \lambda )$ and $\text{dn}(x+\phi | \lambda
)$ become S($x)$ and D($x)$ respectively,~then S,~C,~and D meet the follow equations$^{[8]}$
\begin{eqnarray}
\label{eq5}
\frac{dS}{dx}=CD,~S^2+C^2=1,~D^2+\lambda S^2=1,
\end{eqnarray}
where $\phi $ can determined by equations
\begin{eqnarray}
\label{eq6}
\text{sn}(\phi | \lambda )=S(0),~\text{cn}(\phi | \lambda
)=C(0),~\text{dn}(\phi | \lambda )=D(0).
\end{eqnarray}
If two Jacobi elliptic function solutions of a NDE satisfy  (1) and (2),~we can say that the two solutions are equivalence.~By the means of appendix and  (1) ,~we can get the follow expressions (some of them were given by Ref. ${[8]}$)
\begin{eqnarray}
\label{eq7}
\text{sn}(\sqrt k x| \frac{1}{k})=\sqrt k \text{sn}(x| k),
\end{eqnarray}
\begin{eqnarray}
\label{eq8}
\text{cn}(\sqrt k x| \frac{1}{k})=\text{dn}(x| k),
\end{eqnarray}
\begin{eqnarray}
\label{eq9}
\text{dn}(\sqrt k x| \frac{1}{k})=\text{cn}(x| k).
\end{eqnarray}
\begin{eqnarray}
\label{eq10}
\text{sn}\left( {\sqrt k x+i\sqrt k {K}'(k)+\varphi | \frac{1}{k}}
\right)=(-1)^p\text{ns}(x| k),
\end{eqnarray}
\begin{eqnarray}
\label{eq11}
\text{cn}\left( {\sqrt k x+i\sqrt k {K}'(k)+\varphi | \frac{1}{k}}
\right)=(-1)^{p+q+1}i\text{cs}(x| k),
\end{eqnarray}
\begin{eqnarray}
\label{eq12}
\text{dn}\left( {\sqrt k x+i\sqrt k {K}'(k)+\varphi | \frac{1}{k}}
\right)=\frac{(-1)^{q+1}i}{\sqrt k }\text{ds}(x| k),
\end{eqnarray}
\begin{eqnarray}
 \mbox{sn}\left( {\sqrt k x+\sqrt k K(k)+\varphi | \frac{1}{k}}
\right)=(-1)^p\sqrt k \cdot \mbox{cd}(x| k),~\\
 \mbox{cn}\left( {\sqrt k x+\sqrt k K(k)+\varphi| \frac{1}{k}}
\right)=(-1)^{p+q}\sqrt {{k}'} \cdot \mbox{nd}(x|k),~\\
 \mbox{dn}\left( {\sqrt k x+\sqrt k K(k)+\varphi | \frac{1}{k}}
\right)=(-1)^{q+1}\sqrt {{k}'} \mbox{sd}(x| k),
\end{eqnarray}
\begin{eqnarray}
\label{eq16}
\text{sn}\left( {\sqrt k x+i\sqrt k {K}'(k)+\sqrt k K(k)+\varphi |
\frac{1}{k}} \right)=(-1)^p\cdot \text{dc}(x| k),~\end{eqnarray}
\begin{eqnarray}
\text{cn}\left( {\sqrt k x+i\sqrt k {K}'(k)+\sqrt k K(k)+\varphi |
\frac{1}{k}} \right)=(-1)^{p+q}i\sqrt {{k}'} \cdot \text{cs}(x| k),~\end{eqnarray}
\begin{eqnarray}
\text{dn}\left( {\sqrt k x+i\sqrt k {K}'(k)+\sqrt k K(k)+\varphi |
\frac{1}{k}} \right)=(-1)^{q+1}i\sqrt {\frac{{k}'}{k}} \text{nc}(x| k),
\end{eqnarray}

where $\varphi=2pK(1/k)+2iq{K(1/k)}'$.
\begin{eqnarray}
\label{eq17}
\text{sn}(ix| {k}')=i\cdot \text{sc}(x| k),
\end{eqnarray}
\begin{eqnarray}
\label{eq18}
\text{cn}(ix| {k}')=\text{nc}(x| k),
\end{eqnarray}
\begin{eqnarray}
\label{eq19}
\text{dn}(ix| {k}')=\text{dc}(x| k),
\end{eqnarray}
\begin{eqnarray}
\label{eq20}
\text{sn}(i\sqrt k x| -\frac{{k}'}{k})=i\sqrt k \cdot \text{sd}(x|
k),
\end{eqnarray}
\begin{eqnarray}
\label{eq21}
\text{cn}(i\sqrt k x| -\frac{{k}'}{k})=\text{nd}(x| k),
\end{eqnarray}
\begin{eqnarray}
\label{eq22}
\text{dn}(i\sqrt k x| -\frac{{k}'}{k})=\text{cd}(x| k),
\end{eqnarray}
\begin{eqnarray}
\label{eq23}
\text{sn}(\sqrt {{k}'} x| -\frac{k}{{k}'})=\sqrt {{k}'} \text{sd}(x|
k),\\
\text{cn}(\sqrt {{k}'} x| -\frac{k}{{k}'})=\text{cd}(x| k),
\end{eqnarray}
\begin{eqnarray}
\label{eq24}
\text{dn}(\sqrt {{k}'} x| -\frac{k}{{k}'})=\text{nd}(x| k),
\end{eqnarray}
\begin{eqnarray}
\label{eq25}
\text{sn}(i\sqrt {{k}'} x| \frac{1}{{k}'})=i\sqrt {{k}'}
\text{sc}(x| k),
\end{eqnarray}
\begin{eqnarray}
\label{eq26}
\text{cn}(i\sqrt {{k}'} x| \frac{1}{{k}'})=\text{dc}(x| k),
\end{eqnarray}
\begin{eqnarray}
\label{eq27}
\text{dn}(i\sqrt {{k}'} x| \frac{1}{{k}'})=\text{nc}(x| k),
\end{eqnarray}
\begin{eqnarray}
\label{eq28}
\text{sn}\left( {\left( {\sqrt {{k}'} -\varepsilon } \right)x| \left(
{\frac{\sqrt {{k}'} +\varepsilon }{\sqrt {{k}'} -\varepsilon }} \right)^2}
\right)=\left( {\sqrt {{k}'} -\varepsilon } \right)\text{cd}(x| k)\cdot
\text{sn}(x| k),
\end{eqnarray}
where $\varepsilon^{2}=1$,
\begin{eqnarray}
\label{eq29}
\text{cn}\left( {\left( {\sqrt {{k}'} -\varepsilon } \right)x| \left(
{\frac{\sqrt {{k}'} +\varepsilon }{\sqrt {{k}'} -\varepsilon }} \right)^2}
\right)=\left( {\varepsilon \sqrt {{k}'} -1} \right)\text{sd}(x|
k)\text{sn}(x| k)+\text{nd}(x| k),~\end{eqnarray}
\begin{eqnarray}
\text{dn}\left( {\left( {\sqrt {{k}'} -\varepsilon } \right)x| \left(
{\frac{\sqrt {{k}'} +\varepsilon }{\sqrt {{k}'} -\varepsilon }} \right)^2}
\right)=\frac{1-\left[ {\text{dn}(x| k)^2+1} \right]\text{sn}(x|
k)^2}{\text{dn}(x| k)\left[ {\left( {\varepsilon \sqrt {{k}'} -1}
\right)\text{sn}(x| k))^2+1} \right]},\end{eqnarray}

\begin{eqnarray}
\label{eq34}
\text{sn}\left( {\frac{1}{2}\left( {\varepsilon \sqrt {{k}'} +1}
\right)x+\phi | \left( {\frac{1-\varepsilon \sqrt {{k}'}
}{1+\varepsilon \sqrt {{k}'} }} \right)^2} \right)=\left( {\varepsilon
+\sqrt {{k}'} } \right)\cdot \frac{\text{cn}(x| k)}{\sqrt {{k}'}
-\text{dn}(x| k)},
\end{eqnarray}
\begin{eqnarray}
\label{eq35}
\text{cn}\left( {\frac{1}{2}\left( {\varepsilon \sqrt {{k}'} +1}
\right)x+\phi | \left( {\frac{1-\varepsilon \sqrt {{k}'}
}{1+\varepsilon \sqrt {{k}'} }} \right)^2} \right) 
=i\varepsilon \sqrt
{\frac{2}{k}} \sqrt[4]{{k}'}\frac{\sqrt {\left[ {\text{dn}(x|
k)+\varepsilon } \right]\left\{ {\left( {\sqrt {{k}'} +\varepsilon }
\right)[\text{dn}(x| k)-\varepsilon ]+k} \right\}} }{\text{dn}(x|
k)-\sqrt {{k}'} },
\end{eqnarray}
\begin{eqnarray}
\label{eq36}
\text{dn}\left( {\frac{1}{2}\left( {\varepsilon \sqrt {{k}'} +1}
\right)x+\phi | \left( {\frac{1-\varepsilon \sqrt {{k}'}
}{1+\varepsilon \sqrt {{k}'} }} \right)^2} \right)
=\frac{i\sqrt {2k}
\sqrt[4]{{k}'}\text{sn}(x| k)}{\sqrt {\left[ {\text{dn}(x|
k)+\varepsilon } \right]\left\{ {\left( {\sqrt {{k}'} +\varepsilon }
\right)[\text{dn}(x| k)-\varepsilon ]+k} \right\}} },
\end{eqnarray}
where  $\phi =-\varepsilon K+\frac{1}{2}(1+\varepsilon )i{K}'$,

\begin{eqnarray}
\text{sn}\left( {\frac{i}{2}\left( {\sqrt k +\varepsilon } \right)x+\phi
| \left( {\frac{\sqrt k -\varepsilon }{\sqrt k +\varepsilon }}
\right)^2} \right)=\varepsilon \sqrt {\frac{\varepsilon +\sqrt k
}{\varepsilon -\sqrt k }} \frac{\text{dn}(x| k)}{\sqrt k
\text{sn}(x| k)-1},
 \end{eqnarray}
\begin{eqnarray}\text{cn}\left( {\frac{i}{2}\left( {\sqrt k
+\varepsilon } \right)x+\phi | \left( {\frac{\sqrt k -\varepsilon
}{\sqrt k +\varepsilon }} \right)^2} \right)=\frac{\varepsilon\sqrt {\left(
{\varepsilon +\sqrt k } \right)\text{dn}(x| k)^2+\left( {\sqrt k
-\varepsilon } \right)\left[ {\sqrt k \text{sn}(x| k)-1} \right]^2}
}{\sqrt {\sqrt k -\varepsilon } \left( {1-\sqrt k \text{sn}(x| k)}
\right)},~\quad _{ } \end{eqnarray}
\begin{eqnarray}
\text{dn}\left( {\frac{i}{2}\left( {\sqrt k +\varepsilon } \right)x+\phi
| \left( {\frac{\sqrt k -\varepsilon }{\sqrt k +\varepsilon }}
\right)^2} \right)=\frac{\text{cn}(x| k)}{\varepsilon \sqrt k +1}\sqrt
{\frac{k\left( {\sqrt k +\varepsilon } \right)}{\left( {\varepsilon +\sqrt k
} \right)\text{dn}(x| k)^2+\left( {\sqrt k -\varepsilon } \right)\left[
{\sqrt k \text{sn}(x| k)-1} \right]^2}},
 \end{eqnarray}
where $\phi =\varepsilon \left( {-K+\frac{i{K}'}{2}} \right)$ ,
\begin{eqnarray}
\label{eq41}
\text{sn}\left( {\frac{1}{2}\left( {i\varepsilon \sqrt {{k}'} +\sqrt k }
\right)x|\left( -\sqrt{k} \varepsilon +i \sqrt{1-k}\right)^4}
\right)=\left( {i\varepsilon \sqrt {{k}'} +\sqrt k }
\right)\frac{\text{sn}(x| k)}{1+\text{cn}(x| k)},
\end{eqnarray}
\begin{eqnarray}
\label{eq42}
\text{cn}\left( {\frac{1}{2}\left( {i\varepsilon \sqrt {{k}'} +\sqrt k }
\right)x|\left( -\sqrt{k} \varepsilon +i \sqrt{1-k}\right)^4}
\right)=\sqrt {\frac{2\left( {\sqrt {{k}'} -i\varepsilon \sqrt k }
\right)\left[ {\sqrt {{k}'} +i\varepsilon \sqrt k \text{cn}(x| k)}
\right]}{1+\text{cn}(x| k)}} ,
\end{eqnarray}
\begin{eqnarray}
\text{dn}\left( {\frac{1}{2}\left( {i\varepsilon \sqrt {{k}'} +\sqrt k }
\right)x|\left( -\sqrt{k} \varepsilon +i \sqrt{1-k}\right)^4}
\right)=\frac{\sqrt 2 \text{dn}(x| k)}{\sqrt {\left( {\sqrt {{k}'}
-i\varepsilon \sqrt k } \right)\left[ {\sqrt {{k}'} +i\varepsilon \sqrt k
\text{cn}(x| k)} \right]\left[ {1+\text{cn}(x| k)} \right]} },
\end{eqnarray}
\begin{eqnarray}
\label{eq43}
\text{sn}\left( {2\sqrt {\varepsilon \sqrt k } x+\phi |
\frac{\varepsilon \left( {\varepsilon +\sqrt k } \right)^2}{4\sqrt k }}
\right)=\frac{i}{\left( {\varepsilon +\sqrt k } \right)}\left( {\varepsilon
\sqrt k \text{sn}(x| k)+\text{ns}(x| k)} \right),
\end{eqnarray}
\begin{eqnarray}
\label{eq44}
\text{cn}\left( {2\sqrt {\varepsilon \sqrt k } x+\phi |
\frac{\varepsilon \left( {\varepsilon +\sqrt k } \right)^2}{4\sqrt k }}
\right)=\frac{-i}{\left( {1+\varepsilon \sqrt k } \right)}\left(
{\varepsilon \sqrt k \text{cs}(x| k)\text{dn}(x| k)} \right),
\end{eqnarray}
\begin{eqnarray}
\label{eq45}
\text{dn}\left( {2\sqrt {\varepsilon \sqrt k } x+\phi |
\frac{\varepsilon \left( {\varepsilon +\sqrt k } \right)^2}{4\sqrt k }}
\right)=\frac{\sqrt {-\varepsilon } }{2\sqrt[4]{k}}\left( {\sqrt k
\text{sn}(x| k)-\varepsilon \text{ns}(x| k)} \right),
\end{eqnarray}
where  $\phi =\varepsilon i{K}'$,
\begin{eqnarray}
\text{sn}\left( {2\sqrt \varepsilon
\sqrt[4]{-k{k}'}x+\phi | \frac{\varepsilon \left( {\sqrt {-{k}'}
-\varepsilon \sqrt k } \right)^2}{4\sqrt {-k{k}'} }} \right)
=\left(
{\varepsilon \sqrt {-{k}'} +\sqrt k } \right)\left( {\sqrt k \text{cn}(x|
k)-\varepsilon \sqrt {-{k}'} \text{nc}(x| k)} \right),
~\end{eqnarray}
\begin{eqnarray}
\label{eq46}
\text{cn}\left( {2\sqrt \varepsilon \sqrt[4]{-k{k}'}x+\phi |
\frac{\varepsilon \left( {\sqrt {-{k}'} -\varepsilon \sqrt k }
\right)^2}{4\sqrt {-k{k}'} }} \right)=-\left( {\sqrt k +\varepsilon \sqrt
{-{k}'} } \right)\text{sn}(x| k)\text{dc}(x| k),
\end{eqnarray}
\begin{eqnarray}
\label{eq47}
\text{dn}\left( {2\sqrt \varepsilon \sqrt[4]{-k{k}'}x+\phi |
\frac{\varepsilon \left( {\sqrt {-{k}'} -\varepsilon \sqrt k }
\right)^2}{4\sqrt {-k{k}'} }} \right)=\frac{\varepsilon ^{1/2}\sqrt k
}{2\sqrt[4]{-{k}'k}}\left( {\text{cn}(x| k)+\varepsilon \sqrt
{\frac{-{k}'}{k}} \text{nc}(x| k)} \right),
\end{eqnarray}
where $\phi =\varepsilon K$,
\begin{eqnarray}
\label{eq48}
\text{sn}\left( {2\sqrt {\varepsilon \sqrt {{k}'} } x+\phi|
\frac{-\varepsilon \left( {\sqrt {{k}'} -\varepsilon } \right)^2}{4\sqrt
{{k}'} }} \right)=\frac{-1}{\varepsilon -\sqrt {{k}'} }\left[
{\varepsilon \text{dn}(x| k)-\sqrt {{k}'} \text{nd}(x| k)} \right],
\end{eqnarray}
\begin{eqnarray}
\label{eq49}
\text{cn}\left( {2\sqrt {\varepsilon \sqrt {{k}'} } xx+\phi  |
\frac{-\varepsilon \left( {\sqrt {{k}'} -\varepsilon } \right)^2}{4\sqrt
{{k}'} }} \right)=\left( {\sqrt {{k}'} +\varepsilon }
\right)\text{sn}(x| k)\text{cd}(x| k),
\end{eqnarray}
\begin{eqnarray}
\label{eq50}
\text{dn}\left( {2\sqrt {\varepsilon \sqrt {{k}'} } x+\phi |
\frac{-\varepsilon \left( {\sqrt {{k}'} -\varepsilon } \right)^2}{4\sqrt
{{k}'} }} \right)=\frac{\varepsilon ^{3/2}}{2\sqrt[4]{{k}'}}\left[
{\varepsilon \text{dn}(x| k)+\sqrt {{k}'} \text{nd}(x| k)} \right],
\end{eqnarray}
where  $\phi =-K$,
\begin{eqnarray}
\label{eq51}
\text{sn}\left( {\frac{1}{2}i\varepsilon \left( {\sqrt[4]{k}+\varepsilon
^{3/2}} \right)^2x+\phi | \left( {\frac{\sqrt[4]{k}-\varepsilon
^{3/2}}{\sqrt[4]{k}+\varepsilon ^{3/2}}} \right)^4} \right)=\frac{-\left(
{\sqrt[4]{k}+\varepsilon ^{3/2}} \right)}{\sqrt[4]{k}+\varepsilon
^{3/2}}\frac{\varepsilon ^{3/2}+\sqrt[4]{k}\text{sn}(x| k)}{\varepsilon
^{3/2}-\sqrt[4]{k}\text{sn}(x| k)},
\end{eqnarray}
\begin{eqnarray}
\text{cn}\left( {\frac{1}{2}i\varepsilon \left( {\sqrt[4]{k}+\varepsilon
^{3/2}} \right)^2x+\phi | \left( {\frac{\sqrt[4]{k}-\varepsilon
^{3/2}}{\sqrt[4]{k}+\varepsilon ^{3/2}}} \right)^4}
\right)
=\frac{2i\varepsilon ^{3/4}\sqrt[8]{k}\sqrt {\left(
{1+\text{sn}(x| k)} \right)\left( {\varepsilon +\sqrt k \text{sn}(x|
k)} \right)} }{\left( {\sqrt[4]{k}-\varepsilon ^{3/2}} \right)\left(
{\varepsilon ^{3/2}-\sqrt[4]{k}\text{sn}(x| k)} \right)},
 \end{eqnarray}
\begin{eqnarray}
\text{dn}\left( {\frac{1}{2}i\varepsilon \left( {\sqrt[4]{k}+\varepsilon
^{3/2}} \right)^2x+\phi | \left( {\frac{\sqrt[4]{k}-\varepsilon
^{3/2}}{\sqrt[4]{k}+\varepsilon ^{3/2}}} \right)^4} \right)
= \frac{2\varepsilon ^{-1/4}\sqrt[8]{k}\text{cn}(x|
k)\text{dn}(x| k)}{\left( {\sqrt[4]{k}+\varepsilon ^{3/2}} \right)\left(
{\varepsilon ^{3/2}-\sqrt[4]{k}\text{sn}(x| k)} \right)\sqrt {\left(
{1+\text{sn}(x| k)} \right)\left( {\varepsilon +\sqrt k \text{sn}(x|
k)} \right)} },
 \end{eqnarray}
where  $\phi =\frac{\varepsilon +1}{2}K+\frac{i{K}'}{2}$,
\begin{eqnarray}
\text{sn}\left( {\frac{1}{2}i\left( {\left( {-\varepsilon }
\right)^{3/2}\sqrt[4]{{k}'}+\sqrt[4]{-k}} \right)^2x+\phi | \left(
{\frac{\sqrt[4]{-k}-\left( {-\varepsilon }
\right)^{3/2}\sqrt[4]{{k}'}}{\sqrt[4]{-k}+\left( {-\varepsilon }
\right)^{3/2}\sqrt[4]{{k}'}}} \right)^4} \right)
=-\frac{\left( {\sqrt[4]{-k}+\left( {-\varepsilon }
\right)^{3/2}\sqrt[4]{{k}'}} \right)}{\sqrt[4]{-k}-\left( {-\varepsilon }
\right)^{3/2}\sqrt[4]{{k}'}}\frac{\left( {\left( {-\varepsilon }
\right)^{3/2}\sqrt[4]{{k}'}+\sqrt[4]{-k}\text{cn}(x| k)} \right)}{\left(
{\left( {-\varepsilon } \right)^{3/2}\sqrt[4]{{k}'}-\sqrt[4]{-k}\text{cn}(x|
k)} \right)},
\end{eqnarray}
\begin{eqnarray}
 &&\text{cn}\left( {\frac{1}{2}i\left( {\left( {-\varepsilon }
\right)^{3/2}\sqrt[4]{{k}'}+\sqrt[4]{-k}} \right)^2x+\phi | \left(
{\frac{\sqrt[4]{-k}-\left( {-\varepsilon }
\right)^{3/2}\sqrt[4]{{k}'}}{\sqrt[4]{-k}+\left( {-\varepsilon }
\right)^{3/2}\sqrt[4]{{k}'}}} \right)^4} \right)\nonumber\\
&&~~=\frac{2i\sqrt[8]{-k{k}'}\sqrt {\left(
{2k-\varepsilon } \right)^2} \sqrt {\left( {\text{cn}(x| k)+1}
\right)\left( {\sqrt k \text{cn}(x| k)-\varepsilon \sqrt {{k}'} }
\right)} }{\left( {-\varepsilon } \right)^{1/4}\left( {2k-\varepsilon }
\right)\left( {\sqrt[4]{-k}-\left( {-\varepsilon }
\right)^{3/2}\sqrt[4]{{k}'}} \right)\left( {\sqrt[4]{-k}\text{cn}(x|
k)-\left( {-\varepsilon } \right)^{3/2}\sqrt[4]{{k}'}} \right)},
 \end{eqnarray}
\begin{eqnarray}
 &&\text{dn}\left( {\frac{1}{2}i\left( {\left( {-\varepsilon }
\right)^{3/2}\sqrt[4]{{k}'}+\sqrt[4]{-k}} \right)^2x+\phi | \left(
{\frac{\sqrt[4]{-k}-\left( {-\varepsilon }
\right)^{3/2}\sqrt[4]{{k}'}}{\sqrt[4]{-k}+\left( {-\varepsilon }
\right)^{3/2}\sqrt[4]{{k}'}}} \right)^4} \right)\nonumber\\
&&~~
=\frac{1}{\sqrt {\left(
{\text{cn}(x| k)+1} \right)\left( {\sqrt {-k} \text{cn}(x|
k)-\varepsilon \sqrt {{k}'} } \right)} }
\frac{-2\left( {-\varepsilon }
\right)^{-1/4}\sqrt[8]{-k{k}'}\sqrt {\left( {2k-1} \right)^2}
\text{sn}(x| k)\text{dn}(x| k)}{\left( {2k-\varepsilon }
\right)\left( {\sqrt[4]{-k}+\left( {-\varepsilon }
\right)^{3/2}\sqrt[4]{{k}'}} \right)\left( {\sqrt[4]{-k}\text{cn}(x|
k)-\left( {-\varepsilon } \right)^{3/2}\sqrt[4]{{k}'}} \right)},~\end{eqnarray}
where  $\phi =-K-i{K}'$,
\begin{eqnarray}
\label{eq53}
\text{sn}\left( {\frac{1}{2}\left( {1+\varepsilon ^{3/2}\sqrt[4]{{k}'}}
\right)^2x+\phi | \left( {\frac{1+\varepsilon
^{3/2}\sqrt[4]{{k}'}}{1-\varepsilon ^{3/2}\sqrt[4]{{k}'}}} \right)^4}
\right)=-\varepsilon\frac{\left( {1+\varepsilon ^{3/2}\sqrt[4]{{k}'}} \right)}{\left(
{1-\varepsilon ^{3/2}\sqrt[4]{{k}'}} \right)}\frac{\text{dn}(x|
k)+\varepsilon ^{3/2}\sqrt[4]{{k}'}}{\text{dn}(x| k)-\varepsilon
^{3/2}\sqrt[4]{{k}'}},
\end{eqnarray}
\begin{eqnarray}
\label{eq54}
\text{cn}\left( {\frac{1}{2}\left( {1+\varepsilon ^{3/2}\sqrt[4]{{k}'}}
\right)^2x+\phi | \left( {\frac{1+\varepsilon
^{3/2}\sqrt[4]{{k}'}}{1-\varepsilon ^{3/2}\sqrt[4]{{k}'}}} \right)^4}
\right)=\frac{2i\varepsilon ^{3/4}\sqrt[8]{{k}'}\sqrt {\left(
{1+\text{dn}(x| k)} \right)\left( {\varepsilon \sqrt {{k}'}
+\text{dn}(x| k)} \right)} }{\left( {1-\varepsilon ^{3/2}\sqrt[4]{{k}'}}
\right)\left( {\text{dn}(x| k)-\varepsilon ^{3/2}\sqrt[4]{{k}'}}
\right)},
\end{eqnarray}
\begin{eqnarray}
\label{eq55}
\begin{array}{l}
\text{dn}\left( {\frac{1}{2}\left( {1+\varepsilon ^{3/2}\sqrt[4]{{k}'}}
\right)^2x+\phi | \left( {\frac{1+\varepsilon
^{3/2}\sqrt[4]{{k}'}}{1-\varepsilon ^{3/2}\sqrt[4]{{k}'}}} \right)^4}
\right)=\frac{2ik\varepsilon ^{-1/4}{k}'^{1/8}\text{sn}(x|
k)\text{cn}(x| k)}{\left( {1+\varepsilon ^{3/2}\sqrt[4]{{k}'}}
\right)\left( {\text{dn}(x| k)-\varepsilon ^{3/2}\sqrt[4]{{k}'}}
\right)\sqrt {\left( {1+\text{dn}(x| k)} \right)\left( {\varepsilon
\sqrt {{k}'} +\text{dn}(x| k)} \right)} },
 \end{array}
\end{eqnarray}
where  $\phi =-K-i{K}'$.

\section{investigation of Jacobi elliptic function solutions to mKdV}
In this section,~we will apply our transform to NDE.~We take the mKdV equation for example.
 The mKdV equation reads
\begin{eqnarray}
u_t+6\nu u^2 u_x+u_{\text{xxx}}=0,
\end{eqnarray}
where $\nu^{2}=1$,~$\nu=1$ for positive mKdV(p-mKdV) equation,~and $\nu=-1$ for negative mKdV(n-mKdV) equation.
If we use the relation:
\begin{eqnarray}
u=2\sqrt{\nu} \frac{\partial \tan ^{-1}(\phi )}{\partial x},
\end{eqnarray}
Thus substitution of (58) in (57) gives:
\begin{eqnarray}
\left(1+\phi ^2\right) \left(\phi _t+\phi _{\text{xxx}}\right)+6 \phi _x \left(\phi ^2-\phi \phi _{\text{xx}}\right)=0,
\end{eqnarray}
where $\phi$ is the function of $x$ and $t$.~Ref.~[5,9--11] solved the equation (59),~and got some Jacobi elliptic function solutions.~Fu$^{[11]}$ obtained p-mKdV's more binary Jacobi elliptic function solutions using a systematical way,~and pointed out the conditions for the p-mKdV's binary Jacobi elliptic function solutions existing.~We will discuss the conditions for real solution of both p-mKdV and n-mKdV.  ~As all we known:
\begin{eqnarray}
\tanh ^{-1}(i \phi )=i \tan ^{-1}(\phi ),~i \tanh ^{-1}(\phi)=\tan ^{-1}(i \phi),
\end{eqnarray}
combining (58),~we can find that if $\phi$ is imaginary number,~then the solution of n-mKdV will be real.~Therefor if considering n-mKdV and p-mKdV's real solutions at the same time,~we can treat $\phi^{2}$ as a real.~So when we known a solution of (59) given by Ref.~[5],~we can get a lot of Jacobi elliptic function solutions of n-mKdV and p-mKdV.~Such as
\begin{eqnarray}
\phi_{0}=\frac{i a \sqrt{k}  \sigma }{c}\text{sn}\left(\left.\xi\right|k\right) \text{sn}\left(\left.\eta\right|m\right),
\end{eqnarray}
where
\begin{eqnarray}
\frac{c^4}{a^4}=\frac{k}{m} ,\nonumber\\ b=a \left((k+1) a^2+3 c^2 (m+1)\right) ,\nonumber\\
   \frac{d}{c}=3 (k+1) a^2+c^2 (m+1),
\end{eqnarray}
and
\begin{eqnarray}
\xi =a x+b t+a_0,~\eta =c x+d t+c_0,~\sigma^2=1.
\end{eqnarray}
By the means of  (21), using the transform $\xi \to \sqrt{{k}'}\xi,~\eta \to \sqrt{{m}'}\eta,~k \to \frac{k}{-{k}'},~m \to \frac{m}{-{m}'}$ ,~we can get
\begin{eqnarray}
\phi_{1}=\sigma  \frac{a}{c} \sqrt{k' k}\text{  }\text{sd}(\text{a0}+b t+a x|k) \text{sd}(\text{c0}+d t+c x|m),
  \end{eqnarray}
  where
  \begin{eqnarray}
\frac{c^4}{a^4}=\frac{k' k}{m' m},~\frac{b}{a} (1-2 k) a^2+3 c^2 (1-2 m),~\frac{d}{c}= (3-6 k) a^2+c^2 (1-2 m).
\end{eqnarray}
The other Jacobi elliptic function solutions of mKdV equation are list in Table I.
So these solutions of mKdV equation are equivalent.

\begin{center}
\begin{longtable}{p{0.9cm}p{7.6cm}p{9.2cm}}\caption{\label{tab:test} Jacobi elliptic function solutions of mKdV}\\ 
\toprule 
 \toprule
 \hline
 \hline
 number & $2\sqrt{\nu} \frac{\partial \tan ^{-1}}{\partial x}$& $\frac{c^4}{a^4},~\frac{b}{a},\frac{d}{c}$\\ \midrule 
  \midrule 
  \hline
 
 \endfirsthead {\bf continue~\ref{tab:test}}\\ 
 \toprule
  \toprule
  \hline
   \hline
  number &  $2\sqrt{\nu} \frac{\partial \tan ^{-1}}{\partial x}$& $\frac{c^4}{a^4},\frac{b}{a},\frac{d}{c}$\\
   \midrule
      \midrule
 \hline
  \endhead \endfoot 
  \bottomrule
    \bottomrule
             \hline
     \endlastfoot
 1&$ \frac{i \sigma  a\text{  }\sqrt{k} \text{sc}(\eta |m) \text{sn}(\xi |k)}{c} $&$  \frac{k}{m'},(k+1) a^2+3 c^2 (m-2),3 (k+1) a^2+c^2 (m-2)  $\\
 2&$ \frac{\text{i$\sigma $} a \sqrt{k'} \text{sc}(\eta |m) \text{sc}(\xi |k) }{c} $&$  \frac{k'}{m'},(k-2) a^2+3 c^2 (m-2),3 (k-2) a^2+c^2 (m-2)  $\\
 3&$ \frac{a \sigma  \sqrt{k} \sqrt{-m}\text{  }\text{sd}(\eta |m) \text{sn}(\xi |k)}{c \sqrt{m}} $&$  -\frac{k}{m m'},(k+1) a^2+3 c^2 (1-2 m),3 (k+1) a^2+c^2 (1-2 m)  $\\
 4&$ \frac{\sigma  a\sqrt{k'}m^{5/2} \text{sc}(\xi |k) \text{sd}(\eta |m) }{c (-m)^{5/2}} $&$  -\frac{k'}{m m'},(k-2) a^2+3 c^2 (1-2 m),3 (k-2) a^2+c^2 (1-2 m)  $\\
 5&$ \frac{ \text{$\sigma $a} \sqrt{k k'}\text{   }\text{sd}(\eta |m) \text{sd}(\xi |k) }{c } $&$  \frac{k k'}{m m'},(1-2 k) a^2+3 c^2 (1-2 m),(3-6 k) a^2+c^2 (1-2 m)  $\\
 6&$ \frac{i a \sqrt{k} \sigma  \text{cn}(\eta |m) \text{sn}(\eta |m) \text{sn}(\xi |k)}{c \text{dn}(\eta |m)} $&$  \frac{k}{m^2},a^2 (k+1)-6 c^2 (m-2),3 a^2 (k+1)-2 c^2 (m-2)  $\\
 7&$ \frac{i \sigma  a\text{  }\sqrt{k'} \text{cn}(\eta |m) \text{sc}(\xi |k) \text{sn}(\eta |m)}{c \text{dn}(\eta |m)} $&$  \frac{k'}{m^2},a^2 (k-2)-6 c^2 (m-2),3 a^2 (k-2)-2 c^2 (m-2)  $\\
 8&$ \frac{a \sigma  \sqrt{k k'}\text{cn}(\eta |m) \text{sd}(\xi |k) \text{sn}(\eta |m) }{c \text{dn}(\eta |m)} $&$  -\frac{k k'}{m^2},a^2 (1-2 k)-6 c^2 (m-2),a^2 (3-6 k)-2 c^2 (m-2)  $\\
\hline
 9&$ \frac{i a k \sigma  \text{cn}(\eta |m) \text{cn}(\xi |k) \text{sn}(\eta |m) \text{sn}(\xi |k)}{c \text{dn}(\eta |m) \text{dn}(\xi |k)} $&$  \frac{k^2}{m^2},-2 (k-2) a^2-6 c^2 (m-2),-6 (k-2) a^2-2 c^2 (m-2)  $\\
 10&$ \frac{2 i a \sqrt{k} \sigma  \text{cn}(\eta |m) \text{sn}(\xi |k) \sigma _m}{c\left(\sqrt{m'}-\text{dn}(\eta |m)\right)} $&$  \frac{16 k}{m^2},a^2 (k+1)-\frac{3}{2} c^2 (m-2),3 a^2 (k+1)-\frac{c^2 }{2}(m-2)  $\\
 11&$ \frac{2 i a \sigma \sigma _m \sqrt{k'} \left(\sigma _m+\sqrt{m'}\right) \text{cn}(\eta |m) \text{sc}(\xi |k) }{c\left(\sqrt{m'}-\text{dn}(\eta |m)\right)} $&$  \frac{16 k'}{m^2},\frac{1}{2} \left(2 a^2 (k-2)-3 c^2 (m-2)\right),3 a^2 (k-2)-\frac{c^2 }{2} (m-2)  $\\
 12&$ \frac{2 a \sigma \sigma _m \sqrt{k k'} \left(\sigma _m+\sqrt{m'}\right) \text{cn}(\eta |m) \text{sd}(\xi |k) }{c\left(\sqrt{m'}-\text{dn}(\eta |m)\right)} $&$  -\frac{16 k k'}{m^2},a^2 (1-2 k)-\frac{3}{2} c^2 (m-2),a^2 (3-6 k)-\frac{1}{2} c^2 (m-2)  $\\
 13&$ \frac{2 i a k \sigma  \sigma _m \left(\sigma _m+\sqrt{m'}\right)\text{cn}(\eta |m) \text{cn}(\xi |k) \text{sn}(\xi |k) }{c \text{dn}(\xi |k) \left(\sqrt{m'}-\text{dn}(\eta |m)\right)} $&$  \frac{16 k^2}{m^2},\frac{1}{2} \left(-4 (k-2) a^2-3 c^2 (m-2)\right),-6 (k-2) a^2-\frac{1}{2} c^2 (m-2)  $\\
 14&$ \frac{i a k \sigma  \sigma _k \sigma _m \left(\sigma _m+\sqrt{m'}\right)\text{cn}(\eta |m) \text{cn}(\xi |k) }{c\left(\sqrt{k'}-\text{dn}(\xi |k)\right) \left(\sqrt{m'}-\text{dn}(\eta |m)\right)} $&$  \frac{k^2}{m^2},\frac{1}{2} \left(a^2 (2-k)-3 c^2 (m-2)\right),-\frac{3}{2} (k-2) a^2-\frac{1}{2} c^2  (m-2)  $\\
 15&$ \frac{2 a \sqrt{k} \sigma \sigma _m \text{dn}(\eta |m) \text{sn}(\xi |k) }{c \sqrt{m'}\left(\sqrt{m} \text{sn}(\eta |m)-1\right) } $&$  \frac{16 k}{m'^2},a^2 (k+1)-\frac{3}{2} c^2 (m+1),3 a^2 (k+1)-\frac{c^2}{2} (m+1)  $\\
 16&$ \frac{2 a \sigma \sigma _m \sqrt{k'} \text{dn}(\eta |m) \text{sc}(\xi |k) }{c \sqrt{m'}\left(\sqrt{m} \text{sn}(\eta |m)-1\right) } $&$  \frac{16 k'}{m'^2},a^2 (k-2)-\frac{3}{2} c^2 (m+1),3 a^2 (k-2)-\frac{1}{2} c^2 (m+1)  $\\
 17&$ \frac{2 i a \sigma \sigma _m \sqrt{k k'} \text{dn}(\eta |m) \text{sd}(\xi |k) }{c \sqrt{m'}\left(\sqrt{m} \text{sn}(\eta |m)-1\right) } $&$  -\frac{16 k k'}{m'^2},a^2 (1-2 k)-\frac{3}{2} c^2 (m+1),a^2 (3-6 k)-\frac{1}{2} c^2 (m+1)  $\\

 18&$ \frac{2 a k \sigma \sigma _m \text{cn}(\xi |k) \text{dn}(\eta |m) \text{sn}(\xi |k) }{c\sqrt{m'} \text{dn}(\xi |k) \left(\sqrt{m} \text{sn}(\eta |m)-1\right) } $&$  \frac{16 k^2}{m'^2},-2 (k-2) a^2-\frac{3}{2} c^2 (m+1),-6 (k-2) a^2-\frac{1}{2} c^2 (m+1)  $\\
 19&$ \frac{a k \sigma  \sigma _k \sigma _m\text{cn}(\xi |k) \text{dn}(\eta |m) }{c \sqrt{m'}\left(\sqrt{m} \text{sn}(\eta |m)-1\right) \left(\sqrt{k'}-\text{dn}(\xi |k)\right) } $&$  \frac{k^2}{m'^2},~\frac{1}{2}a^2 (2-k)-\frac{3}{2} c^2 (m+1),-\frac{3}{2} (k-2) a^2-\frac{c^2}{2} (m+1)  $\\
 20&$ \frac{i a \sigma  \sigma _k \sigma _m \sqrt{k'}\text{dn}(\eta |m) \text{dn}(\xi |k) }{c \sqrt{m'}\left(\sqrt{m} \text{sn}(\eta |m)-1\right) \left(\sqrt{k} \text{sn}(\xi |k)-1\right) } $&$  \frac{k'^2}{m'^2},\frac{1}{2}a^2 (-k-1)-\frac{3}{2} c^2 (m+1),-\frac{3}{2} (k+1) a^2-\frac{c^2}{2} (m+1)  $\\
 21&$ \frac{2 a \sigma  \sigma _m \left(\sqrt{m} \sigma _m+i \sqrt{m'}\right) \sqrt{k} \text{sn}(\eta |m) \text{sn}(\xi |k) }{c (\text{cn}(\eta |m)+1) \left(\sigma _m \sqrt{m'}-i \sqrt{m}\right)} $&$  16 k,(k+1) a^2+\frac{3}{2} c^2 (2 m-1),3 (k+1) a^2+\frac{1}{2} c^2 (2 m-1)  $\\
 22&$ \frac{2 a \sigma  \sigma _m \sqrt{k'} \left(\sqrt{m} \sigma _m+i \sqrt{m'}\right) \text{sc}(\xi |k) \text{sn}(\eta |m)}{c \left(\sigma _m \sqrt{m'}-i \sqrt{m}\right)(\text{cn}(\eta |m)+1) } $&$  16 k',(k-2) a^2+\frac{3}{2} c^2 (2 m-1),3 (k-2) a^2+\frac{1}{2} c^2 (2 m-1)  $\\
 23&$ \frac{2 a \sigma \sigma _m \sqrt{k k'} \left(\sqrt{m} \sigma _m+i \sqrt{m'}\right) \text{sd}(\xi |k) \text{sn}(\eta |m) }{c\left(i \sqrt{m'} \sigma _m+\sqrt{m}\right) (\text{cn}(\eta |m)+1) } $&$  -16 k k',(1-2 k) a^2+\frac{3}{2} c^2 (2 m-1),(3-6 k) a^2+\frac{1}{2} c^2 (2 m-1)  $\\
 24&$ \frac{2 a k \sigma \sigma _m \left(\sqrt{m} \sigma _m+i \sqrt{m'}\right) \text{cn}(\xi |k) \text{sn}(\eta |m) \text{sn}(\xi |k) }{c\left(\sigma _m \sqrt{m'}-i \sqrt{m}\right) (\text{cn}(\eta |m)+1) \text{dn}(\xi |k) } $&$  16 k^2,\frac{3}{2} c^2 (2 m-1)-2 a^2 (k-2),\frac{1}{2} c^2 (2 m-1)-6 a^2 (k-2)  $\\
 25&$ \frac{a k \sigma  \sigma _k \sigma _m \left(\sqrt{m} \sigma _m+i \sqrt{m'}\right)\text{cn}(\xi |k) \text{sn}(\eta |m) }{c\left(\sigma _m \sqrt{m'}-i \sqrt{m}\right) (\text{cn}(\eta |m)+1) \left(\sqrt{k'}-\text{dn}(\xi |k)\right) } $&$  k^2,\frac{3}{2} c^2 (2 m-1)-\frac{1}{2} a^2 (k-2),~\frac{1}{2}c^2 (2 m-1)-\frac{3}{2} a^2 (k-2)  $\\
 26&$ \frac{a \sigma \sigma _k \sigma _m \sqrt{k'} \left(\sqrt{m} \sigma _m+i \sqrt{m'}\right) \text{dn}(\xi |k) \text{sn}(\eta |m) }{c\left(i \sqrt{m'} \sigma _m+\sqrt{m}\right)(\text{cn}(\eta |m)+1) \left(\sqrt{k} \text{sn}(\xi |k)-1\right) } $&$  k'^2,~\frac{3}{2} c^2 (2 m-1)-\frac{1}{2}a^2 (k+1),\frac{1}{2} c^2 (2 m-1)-\frac{3}{2}  a^2 (k+1)  $\\
 27&$ \frac{i a \sigma \sigma _m (1-2 \mathrm{\theta} [0.5-k]) \left(\sqrt{m} \sigma _m+i \sqrt{m'}\right) \text{sn}(\eta |m) \text{sn}(\xi |k) }{c\left(i \sqrt{m'} \sigma _m+\sqrt{m}\right) (\text{cn}(\eta |m)+1) (\text{cn}(\xi |k)+1) } $&$  1,\left(k-\frac{1}{2}\right) a^2+\frac{3}{2} c^2 (2 m-1),\left(3 k-\frac{3}{2}\right) a^2+\frac{1}{2} c^2 (2 m-1)  $\\
 28&$ \frac{i a \sqrt{k} \sigma  \text{sn}(\xi |k) \left(\sqrt{m} \sigma _m \text{sn}(\eta |m)^2+1\right)}{2 c\sqrt{\sigma _m} \left(\sigma _m+\sqrt{m}\right) \sqrt[4]{m} \text{sn}(\eta |m) } $&$  \frac{k}{4 \sqrt{m} \sigma _m \left(\sigma _m+\sqrt{m}\right){}^2},(k+1) a^2+3 c^2 (m+1)+18 c^2 \sqrt{m} \sigma _m,3 (k+1) a^2+c^2 (m+1)+6 c^2 \sqrt{m} \sigma _m  $\\
 29&$\frac{i a \sigma \sqrt{k'} \text{sc}(\xi |k) \left(\sqrt{m} \sigma _m \text{sn}(\eta |m)^2+1\right) }{2 c \sqrt[4]{m}\sqrt{\sigma _m} \left(\sigma _m+\sqrt{m}\right) \text{sn}(\eta |m) } $&$  \frac{k'}{4 \sqrt{m} \sigma _m \left(\sigma _m+\sqrt{m}\right){}^2},(k-2) a^2+3 c^2 (m+1)+18 c^2 \sqrt{m} \sigma _m,3 (k-2) a^2+c^2 (m+1)+6 c^2 \sqrt{m} \sigma _m  $\\
 30&$ \frac{a \sigma \sqrt{k k'} \text{sd}(\xi |k) \left(\sqrt{m} \sigma _m \text{sn}(\eta |m)^2+1\right) }{2 c \sqrt{\sigma _m} \left(\sigma _m+\sqrt{m}\right)\sqrt[4]{m} \text{sn}(\eta |m) } $&$  -\frac{k k'}{4 \sqrt{m} \sigma _m \left(\sigma _m+\sqrt{m}\right){}^2},(1-2 k) a^2+3 c^2 (m+1)+18 c^2 \sqrt{m} \sigma _m,(3-6 k) a^2+c^2 (m+1)+6 c^2 \sqrt{m} \sigma _m  $\\
 31&$ \frac{i a k \sigma  \text{cn}(\xi |k) \text{sn}(\xi |k) \left(\sqrt{m} \sigma _m \text{sn}(\eta |m)^2+1\right)}{2 c\sqrt{\sigma _m} \left(\sigma _m+\sqrt{m}\right) \sqrt[4]{m} \text{dn}(\xi |k) \text{sn}(\eta |m) } $&$  \frac{k^2}{4 \sqrt{m} \sigma _m \left(\sigma _m+\sqrt{m}\right){}^2},-2 (k-2) a^2+3 c^2 (m+1)+18 c^2 \sqrt{m} \sigma _m,-6 (k-2) a^2+c^2 (m+1)+6 c^2 \sqrt{m} \sigma _m  $\\
 32&$ \frac{i a k \sigma  \text{cn}(\xi |k) \sigma _k \left(\sqrt{m} \sigma _m \text{sn}(\eta |m)^2+1\right)}{4 c \sqrt[4]{m}\text{  }\sqrt{\sigma _m} \left(\sigma _m+\sqrt{m}\right)\text{sn}(\eta |m) \left(\sqrt{k'}-\text{dn}(\xi |k)\right)} $&$  \frac{k^2}{64 \sqrt{m} \sigma _m \left(\sigma _m+\sqrt{m}\right){}^2},\frac{1}{2}(2-k) a^2+3 c^2 (m+1)+18 c^2 \sqrt{m} \sigma _m,-\frac{3}{2} (k-2) a^2+c^2 (m+1)+6 c^2 \sqrt{m} \sigma _m  $\\
 33&$ \frac{a \sigma  \sigma _k\sqrt{k'}\text{dn}(\xi |k)\text{  }\left(\sqrt{m} \sigma _m \text{sn}(\eta |m)^2+1\right) }{4 c \sqrt[4]{m}\text{  }\sqrt{\sigma _m} \left(\sigma _m+\sqrt{m}\right)\text{sn}(\eta |m) \left(\sqrt{k} \text{sn}(\xi |k)-1\right)} $&$  \frac{k'^2}{64 \sqrt{m} \sigma _m \left(\sigma _m+\sqrt{m}\right){}^2},~\frac{1}{2}(-k-1) a^2+3 c^2 (m+1)+18 c^2 \sqrt{m} \sigma _m,-\frac{3}{2} (k+1) a^2+c^2 (m+1)+6 c^2 \sqrt{m} \sigma _m  $\\
 34&$ \frac{i a \sigma \left(1-2 \mathrm{\theta} [0.5-k] \text{sn}(\xi |k) \left(\sqrt{m} \sigma _m \text{sn}(\eta |m)^2+1\right) \right)}{4 c \sqrt[4]{m}\sqrt{\sigma _m} \left(\sigma _m+\sqrt{m}\right) (\text{cn}(\xi |k)+1) \text{sn}(\eta |m) } $&$  \frac{1}{64 \sqrt{m} \sigma _m \left(\sigma _m+\sqrt{m}\right){}^2},\left(k-\frac{1}{2}\right) a^2+3 c^2 (m+1)+18 c^2 \sqrt{m} \sigma _m,\left(3 k-\frac{3}{2}\right) a^2+c^2 (m+1)+6 c^2 \sqrt{m} \sigma _m  $\\
 35&$ \frac{i a \sigma  \left(\sqrt{k} \sigma _k \text{sn}(\xi |k)^2+1\right) \left(\sqrt{m} \sigma _m \text{sn}(\eta |m)^2+1\right)}{2 c \sqrt[4]{m} \sqrt{\sigma _m} \left(\sigma _m+\sqrt{m}\right)\text{sn}(\eta |m) \text{sn}(\xi |k) } $&$  \frac{\sqrt{\frac{k}{m}} \sigma _k \left(\sigma _k+\sqrt{k}\right){}^2}{\sigma _m \left(\sigma _m+\sqrt{m}\right){}^2},(k+1) a^2+6 \sqrt{k} \sigma _k a^2+3 c^2 (m+1)+18 c^2 \sqrt{m} \sigma _m,3 (k+1) a^2+18 \sqrt{k} \sigma _k a^2+c^2 (m+1)+6 c^2 \sqrt{m} \sigma _m  $\\
 36&$ \frac{i a \sigma \sqrt{k}\text{   }\sqrt{\sigma _m}\text{  }\left(m+\sigma _m \sqrt{-m m'}\right)}{2 c \sqrt[4]{-m m'}} \text{sn}(\xi |k)\newline\times\left(\text{cn}(\eta |m)-\frac{\sigma _m \sqrt{-\frac{m'}{m}}}{\text{cn}(\eta |m)}\right) $&$  -\frac{k}{4 \sigma _m \left(\sqrt{-m'}-\sqrt{m} \sigma _m\right){}^2 \sqrt{-m m'}},(k+1) a^2+3 c^2 (1-2 m)+18 c^2 \sigma _m \sqrt{-m m'},3 (k+1) a^2+c^2 (1-2 m)+6 c^2 \sigma _m \sqrt{-m m'}  $\\
 37&$ \frac{i a \sigma \text{  }\sqrt{\sigma _m} \sqrt{k'}\left(m+\sigma _m \sqrt{-m m'}\right)}{2 c \sqrt[4]{-m m'}}\text{sc}(\xi |k)\newline\times  \left(\text{cn}(\eta |m)-\frac{\sigma _m \sqrt{-\frac{m'}{m}}}{\text{cn}(\eta |m)}\right)  $&$  -\frac{k'}{4 \sigma _m \left(\sqrt{-m'}-\sqrt{m} \sigma _m\right){}^2 \sqrt{-m m'}},(k-2) a^2+3 c^2 (1-2 m)+18 c^2 \sigma _m \sqrt{-m m'},3 (k-2) a^2+c^2 (1-2 m)+6 c^2 \sigma _m \sqrt{-m m'}  $\\

\hline
 38&$\frac{a \sigma \text{  }\sqrt{\sigma _m} \sqrt{k k'}\text{  }\left(m+\sigma _m \sqrt{-m m'}\right)}{2 c\text{  }\sqrt[4]{-m m'}}\text{sd}(\xi |k)\newline\times \frac{\left(\text{cn}(\eta |m)^2-\sigma _m \sqrt{-\frac{m'}{m}}\right)}{\text{cn}(\eta |m)} $&$  \frac{k k'}{4 \sigma _m \left(\sqrt{-m'}-\sqrt{m} \sigma _m\right){}^2 \sqrt{-m m'}},(1-2 k) a^2+3 c^2 (1-2 m)+18 c^2 \sigma _m \sqrt{-m m'},(3-6 k) a^2+c^2 (1-2 m)+6 c^2 \sigma _m \sqrt{-m m'}  $\\
 39&$ \frac{i a k \sigma  \sqrt{\sigma _m} \left(m+\sigma _m \sqrt{-m m'}\right)\text{cn}(\xi |k) \text{sn}(\xi |k) }{2 c\text{  }\sqrt[4]{-m m'}\text{dn}(\xi |k)}\newline\times \left(\text{cn}(\eta |m)-\frac{\sigma _m \sqrt{-\frac{m'}{m}}}{\text{cn}(\eta |m)}\right)  $&$  -\frac{k^2}{4 \sigma _m \left(\sqrt{-m'}-\sqrt{m} \sigma _m\right){}^2 \sqrt{-m m'}},-2 (k-2) a^2+3 c^2 (1-2 m)+18 c^2 \sigma _m \sqrt{-m m'},-6 (k-2) a^2+c^2 (1-2 m)+6 c^2 \sigma _m \sqrt{-m m'}  $\\
 40&$ \frac{i a k \sigma \text{  }\sigma _k \sqrt{\sigma _m}\text{  }\left(m+\sigma _m \sqrt{-m m'}\right)\text{cn}(\xi |k)}{4 c\text{  }\sqrt[4]{-m m'}\left(\sqrt{k'}-\text{dn}(\xi |k)\right)}\newline\times \left(\text{cn}(\eta |m)-\frac{\sigma _m \sqrt{-\frac{m'}{m}}}{\text{cn}(\eta |m)}\right) $&$  -\frac{k^2}{64 \sigma _m \left(\sqrt{-m'}-\sqrt{m} \sigma _m\right){}^2 \sqrt{-m m'}},-\frac{1}{2} (k-2) a^2+3 c^2 (1-2 m)+18 c^2 \sigma _m \sqrt{-m m'},-\frac{3}{2} (k-2) a^2+c^2 (1-2 m)+6 c^2 \sigma _m \sqrt{-m m'}  $\\
 41&$ \frac{a \sigma \text{  }\sigma _k \sqrt{\sigma _m} \sqrt{k'}\text{  }\left(m+\sigma _m \sqrt{-m m'}\right)\text{dn}(\xi |k)}{4 c\sqrt[4]{-m m'} \left(\sqrt{k} \text{sn}(\xi |k)-1\right) }\newline\times \left(\text{cn}(\eta |m)-\frac{\sigma _m \sqrt{-\frac{m'}{m}}}{\text{cn}(\eta |m)}\right) $&$  -\frac{k'^2}{64 \sigma _m^3 \left(\sqrt{-m'}-\sqrt{m} \sigma _m\right){}^2 \sqrt{-m m'}},-\frac{1}{2} (k+1) a^2+3 c^2 (1-2 m)+18 c^2 \sigma _m \sqrt{-m m'},-\frac{3}{2} (k+1) a^2+c^2 (1-2 m)+6 c^2 \sigma _m \sqrt{-m m'}  $\\
 42&$ \frac{i a \sigma \text{  }\sqrt{\sigma _m} (1-2 \mathrm{\theta} [0.5 -k])\text{  }\left(m+\sigma _m \sqrt{-m m'}\right)\text{sn}(\xi |k)}{4 c \sqrt[4]{-m m'} (\text{cn}(\xi |k)+1)}\newline\times \left(\text{cn}(\eta |m)-\frac{\sigma _m \sqrt{-\frac{m'}{m}}}{\text{cn}(\eta |m)}\right) $&$  -\frac{1}{64 \sigma _m \left(\sqrt{-m'}-\sqrt{m} \sigma _m\right){}^2 \sqrt{-m m'}},\left(k-\frac{1}{2}\right) a^2+3 c^2 (1-2 m)+18 c^2 \sigma _m \sqrt{-m m'},\left(3 k-\frac{3}{2}\right) a^2+c^2 (1-2 m)+6 c^2 \sigma _m \sqrt{-m m'}  $\\
 43&$ \frac{i a \sigma \text{  }\sqrt{\sigma _m}\text{  }\left(m+\sigma _m \sqrt{-m m'}\right)\left(\sqrt{k} \sigma _k \text{sn}(\xi |k)^2+1\right)}{2 c\text{  }\sqrt[4]{-m m'}\text{sn}(\xi |k)}\newline\times \left(\text{cn}(\eta |m)-\frac{\sigma _m \sqrt{-\frac{m'}{m}}}{\text{cn}(\eta |m)}\right) $&$  -\frac{k \sigma _k \left(\sigma _k+\sqrt{k}\right){}^2}{\sigma _m \left(\sqrt{-m'}-\sqrt{m} \sigma _m\right){}^2 \sqrt{-k m m'}},(k+1) a^2+6 \sqrt{k} \sigma _k a^2+3 c^2 (1-2 m)+18 c^2 \sigma _m \sqrt{-m m'},3 (k+1) a^2+18 \sqrt{k} \sigma _k a^2+c^2 (1-2 m)+6 c^2 \sigma _m \sqrt{-m m'}  $\\
 44&$ \frac{i a\text{  }\sigma  \sigma _k \sqrt{k}\sqrt{\sigma _m} \left(m+\sigma _m \sqrt{-m m'}\right)}{2 c \sqrt[4]{-m m'}}\newline\times  \left(\text{cn}(\xi |k)-\frac{\sigma _k \sqrt{-\frac{k'}{k}}}{\text{cn}(\xi |k)}\right) \left(\text{cn}(\eta |m)-\frac{\sigma _m \sqrt{-\frac{m'}{m}}}{\text{cn}(\eta |m)}\right) $&$  \frac{\sigma _k \left(\sqrt{-k'}-\sqrt{k} \sigma _k\right){}^2 \sqrt{-k k'}}{\sigma _m \left(\sqrt{-m'}-\sqrt{m} \sigma _m\right){}^2 \sqrt{-m m'}},(1-2 k) a^2+6 \sigma _k \sqrt{-k k'} a^2+3 c^2 (1-2 m)+18 c^2 \sigma _m \sqrt{-m m'},(3-6 k) a^2+18 \sigma _k \sqrt{-k k'} a^2+c^2 (1-2 m)+6 c^2 \sigma _m \sqrt{-m m'}  $\\
 45&$ \frac{i a \sigma \sqrt{k}\text{   }\sqrt{\sigma _m}\text{sn}(\xi |k) \left(\sqrt{m'}-\text{dn}(\eta |m)^2 \sigma _m\right)}{2 c\text{  }\left(\sigma _m \sqrt{m'}-1\right) \sqrt[4]{m'}\text{dn}(\eta |m)} $&$  -\frac{k}{4 \sigma _m \left(\sqrt{m'}-\sigma _m\right){}^2 \sqrt{m'}},(k+1) a^2+3 c^2 (m-2)+18 c^2 \sigma _m \sqrt{m'},3 (k+1) a^2+c^2 (m-2)+6 c^2 \sigma _m \sqrt{m'}  $\\
 46&$ \frac{i a \sigma \text{  }\sqrt{\sigma _m} \sqrt{k'} \text{sc}(\xi |k)\left(\sqrt{m'}-\text{dn}(\eta |m)^2 \sigma _m\right)}{2 c\text{  }\left(\sigma _m \sqrt{m'}-1\right) \sqrt[4]{m'}\text{dn}(\eta |m)} $&$  -\frac{k'}{4 \sigma _m \left(\sqrt{m'}-\sigma _m\right){}^2 \sqrt{m'}},(k-2) a^2+3 c^2 (m-2)+18 c^2 \sigma _m \sqrt{m'},3 (k-2) a^2+c^2 (m-2)+6 c^2 \sigma _m \sqrt{m'}  $\\
 47&$ \frac{a \sigma \text{  }\sqrt{\sigma _m} \sqrt{k k'}\text{sd}(\xi |k) \left(\sqrt{m'}-\text{dn}(\eta |m)^2 \sigma _m\right)}{2 c\text{  }\left(\sigma _m \sqrt{m'}-1\right) \sqrt[4]{m'}\text{dn}(\eta |m)} $&$  \frac{k k'}{4 \sigma _m \left(\sqrt{m'}-\sigma _m\right){}^2 \sqrt{m'}},(1-2 k) a^2+3 c^2 (m-2)+18 c^2 \sigma _m \sqrt{m'},(3-6 k) a^2+c^2 (m-2)+6 c^2 \sigma _m \sqrt{m'}  $\\
 48&$ \frac{i a k \sigma  \sqrt{\sigma _m}\text{cn}(\xi |k) \text{sn}(\xi |k)\text{  }\left(\sqrt{m'}-\text{dn}(\eta |m)^2 \sigma _m\right)}{2 c \left(\sigma _m \sqrt{m'}-1\right) \sqrt[4]{m'} \text{dn}(\eta |m) \text{dn}(\xi |k)} $&$  -\frac{k^2}{4 \sigma _m \left(\sqrt{m'}-\sigma _m\right){}^2 \sqrt{m'}},-2 (k-2) a^2+3 c^2 (m-2)+18 c^2 \sigma _m \sqrt{m'},-6 (k-2) a^2+c^2 (m-2)+6 c^2 \sigma _m \sqrt{m'}  $\\
 49&$ \frac{i a k \sigma \text{  }\sigma _k \sqrt{\sigma _m}\text{cn}(\xi |k) \left(\text{dn}(\eta |m)^2 \sigma _m-\sqrt{m'}\right)}{4 c\left(\sigma _m \sqrt{m'}-1\right) \sqrt[4]{m'} \text{dn}(\eta |m) \left(\sqrt{k'}-\text{dn}(\xi |k)\right) } $&$  -\frac{k^2}{64 \sigma _m \left(\sqrt{m'}-\sigma _m\right){}^2 \sqrt{m'}},-\frac{1}{2} (k-2) a^2+3 c^2 (m-2)+18 c^2 \sigma _m \sqrt{m'},-\frac{3}{2} (k-2) a^2+c^2 (m-2)+6 c^2 \sigma _m \sqrt{m'}  $\\
 50&$ \frac{a \sigma \text{  }\sigma _k \sqrt{\sigma _m} \sqrt{k'}\text{dn}(\xi |k) \left(\text{dn}(\eta |m)^2 \sigma _m-\sqrt{m'}\right)}{4 c\left(\sigma _m \sqrt{m'}-1\right) \sqrt[4]{m'} \text{dn}(\eta |m) \left(\sqrt{k} \text{sn}(\xi |k)-1\right) } $&$  -\frac{k'^2}{64 \sigma _m \left(\sqrt{m'}-\sigma _m\right){}^2 \sqrt{m'}},-\frac{1}{2} (k+1) a^2+3 c^2 (m-2)+18 c^2 \sigma _m \sqrt{m'},-\frac{3}{2} (k+1) a^2+c^2 (m-2)+6 c^2 \sigma _m \sqrt{m'}  $\\
 51&$ \frac{a \sigma  i(1-2 \mathrm{\theta}[0.5 -k])\text{  }\sqrt{\sigma _m}\text{sn}(\xi |k) \left(\sqrt{m'}-\text{dn}(\eta |m)^2 \sigma _m\right)}{4 c\left(\sigma _m \sqrt{m'}-1\right) \sqrt[4]{m'} (\text{cn}(\xi |k)+1) \text{dn}(\eta |m) } $&$  -\frac{1}{64 \sigma _m \left(\sqrt{m'}-\sigma _m\right){}^2 \sqrt{m'}},\left(k-\frac{1}{2}\right) a^2+3 c^2 (m-2)+18 c^2 \sigma _m \sqrt{m'},\left(3 k-\frac{3}{2}\right) a^2+c^2 (m-2)+6 c^2 \sigma _m \sqrt{m'}  $\\
 52&$ \frac{i a \sigma  \sqrt{\sigma _m}\left(\sqrt{k} \sigma _k \text{sn}(\xi |k)^2+1\right)\text{  }\left(\sqrt{m'}-\text{dn}(\eta |m)^2 \sigma _m\right)}{2 c\left(\sigma _m \sqrt{m'}-1\right) \sqrt[4]{m'} \text{dn}(\eta |m) \text{sn}(\xi |k) } $&$  -\frac{\sqrt{k-k m} \sigma _k \left(\sigma _k+\sqrt{k}\right){}^2}{\sigma _m \left(\sqrt{m'}-\sigma _m\right){}^2 m'},(k+1) a^2+6 \sqrt{k} \sigma _k a^2+3 c^2 (m-2)+18 c^2 \sigma _m \sqrt{m'},3 (k+1) a^2+18 \sqrt{k} \sigma _k a^2+c^2 (m-2)+6 c^2 \sigma _m \sqrt{m'}  $\\
 53&$ \frac{a \sqrt{k} \sigma  \sigma _k \sqrt{\sigma _m} \left(\text{dn}(\eta |m)^2 \sigma _m-\sqrt{m'}\right)}{2 c\text{  }\left(\sigma _m \sqrt{m'}-1\right) \sqrt[4]{m'}\text{dn}(\eta |m)}\newline\times \left(\text{cn}(\xi |k)-\frac{\sigma _k \sqrt{-\frac{k'}{k}}}{\text{cn}(\xi |k)}\right)  $&$  \frac{\sigma _k \left(\sqrt{-k'}-\sqrt{k} \sigma _k\right){}^2 \sqrt{-\frac{k k'}{m'}}}{\sigma _m \left(\sqrt{m'}-\sigma _m\right){}^2},(1-2 k) a^2+6 \sigma _k \sqrt{-k k'} a^2+3 c^2 (m-2)+18 c^2 \sigma _m \sqrt{m'},(3-6 k) a^2+18 \sigma _k \sqrt{-k k'} a^2+c^2 (m-2)+6 c^2 \sigma _m \sqrt{m'}  $\\
 54&$ \frac{a \sigma  \sigma _k \sqrt{\sigma _m} \left(\sqrt{k'}-\text{dn}(\xi |k)^2 \sigma _k\right) \left(\sqrt{m'}-\text{dn}(\eta |m)^2 \sigma _m\right)}{2 c \left(\sigma _m \sqrt{m'}-1\right) \sqrt[4]{m'} \text{dn}(\eta |m) \text{dn}(\xi |k)} $&$  \frac{\sigma _k \left(\sqrt{k'}-\sigma _k\right){}^2 \sqrt{\frac{k'}{m'}}}{\sigma _m \left(\sqrt{m'}-\sigma _m\right){}^2},(k-2) a^2+6 \sigma _k \sqrt{k'} a^2+3 c^2 (m-2)+18 c^2 \sigma _m \sqrt{m'},3 (k-2) a^2+18 \sigma _k \sqrt{k'} a^2+c^2 (m-2)+6 c^2 \sigma _m \sqrt{m'}  $\\
 55&$ \frac{2 a \sqrt{k} \sigma  \text{sn}(\xi |k) \left(\sigma _m^{3/2}+\sqrt[4]{m} \text{sn}(\eta |m)\right)}{c \left(\sqrt{m}-\sigma _m\right) \sigma _m \left(\sigma _m^{3/2}-\sqrt[4]{m} \text{sn}(\eta |m)\right)} $&$\frac{16 k}{\left(\sqrt{m} \sigma _m-1\right){}^4},(k+1) a^2-\frac{3}{2} c^2 (m+1)-9 c^2 \sqrt{m} \sigma _m,3 (k+1) a^2-\frac{c^2}{2} (m+1)-3c^2 \sqrt{m} \sigma _m$\\

\hline
 56&$ \frac{2 a \sigma \sqrt{k'} \text{sc}(\xi |k) \left(\sigma _m^{3/2}+\sqrt[4]{m} \text{sn}(\eta |m)\right) }{c \left(\sqrt{m}-\sigma _m\right) \sigma _m \left(\sigma _m^{3/2}-\sqrt[4]{m} \text{sn}(\eta |m)\right)} $&$  \frac{16 k'}{\left(\sqrt{m} \sigma _m-1\right){}^4},(k-2) a^2-\frac{3}{2} c^2 (m+1)-9 c^2 \sqrt{m} \sigma _m,3 (k-2) a^2-\frac{1}{2} c^2 (m+1)-3 c^2 \sqrt{m} \sigma _m  $\\
 57&$ \frac{2 i a \sigma  \sqrt{k k'}\text{sd}(\xi |k) \left(\sigma _m^{3/2}+\sqrt[4]{m} \text{sn}(\eta |m)\right) }{c \left(\sqrt{m}-\sigma _m\right) \sigma _m \left(\sigma _m^{3/2}-\sqrt[4]{m} \text{sn}(\eta |m)\right)} $&$  -\frac{16 k k'}{\left(\sqrt{m} \sigma _m-1\right){}^4},(1-2 k) a^2-\frac{3}{2} c^2 (m+1)-9 c^2 \sqrt{m} \sigma _m,(3-6 k) a^2-\frac{1}{2} c^2 (m+1)-3 c^2 \sqrt{m} \sigma _m  $\\
 58&$ \frac{2 a k \sigma  \text{cn}(\xi |k) \text{sn}(\xi |k) \left(\sigma _m^{3/2}+\sqrt[4]{m} \text{sn}(\eta |m)\right)}{c\text{  }\left(\sqrt{m}-\sigma _m\right) \sigma _m \text{dn}(\xi |k)\left(\sigma _m^{3/2}-\sqrt[4]{m} \text{sn}(\eta |m)\right)} $&$  \frac{16 k^2}{\left(\sqrt{m} \sigma _m-1\right){}^4},~\left(-2 (k-2) a^2-\frac{3}{2} c^2 (m+1)\right)-9 c^2 \sqrt{m} \sigma _m,-6 (k-2) a^2-\frac{1}{2} c^2 (m+1)-3 c^2 \sqrt{m} \sigma _m  $\\
 59&$ \frac{a k \sigma \sigma _k \text{cn}(\xi |k)\text{  }\left(\sigma _m^{3/2}+\sqrt[4]{m} \text{sn}(\eta |m)\right)}{c \left(\sqrt{m}-\sigma _m\right) \sigma _m \left(\sigma _m^{3/2}-\sqrt[4]{m} \text{sn}(\eta |m)\right) \left(\sqrt{k'}-\text{dn}(\xi |k)\right)} $&$  \frac{k^2}{\left(\sqrt{m} \sigma _m-1\right){}^4},\frac{1}{2}(2-k) a^2-\frac{3}{2} c^2 (m+1)-9c^2 \sqrt{m} \sigma _m,-\frac{3}{2} (k-2) a^2-\frac{c^2 }{2}(m+1)-3c^2 \sqrt{m} \sigma _m  $\\
 60&$ \frac{i a \sigma \sigma _k \sqrt{k'} \text{dn}(\xi |k) \left(\sigma _m^{3/2}+\sqrt[4]{m} \text{sn}(\eta |m)\right) }{c\text{  }\left(\sqrt{m}-\sigma _m\right) \sigma _m \left(\sqrt{k} \text{sn}(\xi |k)-1\right)\left(\sigma _m^{3/2}-\sqrt[4]{m} \text{sn}(\eta |m)\right)} $&$  \frac{k'^2}{\left(\sqrt{m} \sigma _m-1\right){}^4},~\frac{1}{2}(-k-1) a^2-\frac{3}{2} c^2 (m+1)-9c^2 \sqrt{m} \sigma _m,-\frac{3}{2} (k+1) a^2-\frac{c^2 }{2}(m+1)-3 c^2 \sqrt{m} \sigma _m  $\\
 61&$ \frac{a \sigma  (1-2 \mathrm{\theta} [0.5 -k]) \text{sn}(\xi |k) \left(\sigma _m^{3/2}+\sqrt[4]{m} \text{sn}(\eta |m)\right)}{c\left(\sqrt{m}-\sigma _m\right) \sigma _m (\text{cn}(\xi |k)+1)\text{  }\left(\sigma _m^{3/2}-\sqrt[4]{m} \text{sn}(\eta |m)\right)} $&$  \frac{1}{\left(\sqrt{m} \sigma _m-1\right){}^4},\frac{1}{2}(2 k-1) a^2-\frac{3}{2} c^2 (m+1)-9 c^2 \sqrt{m} \sigma _m,\frac{1}{2}(6 k-3) a^2-\frac{c^2}{2} (m+1)-3 c^2 \sqrt{m} \sigma _m  $\\
 62&$ \frac{2 a \sigma  \left(\sqrt{k} \sigma _k \text{sn}(\xi |k)^2+1\right) \left(\sigma _m^{3/2}+\sqrt[4]{m} \text{sn}(\eta |m)\right)}{c \sigma _m \left(\sqrt{m}-\sigma _m\right)\text{  }\text{sn}(\xi |k)\left(\sigma _m^{3/2}-\sqrt[4]{m} \text{sn}(\eta |m)\right)} $&$  \frac{64 \sqrt{k} \sigma _k \left(\sigma _k+\sqrt{k}\right){}^2}{\left(\sqrt{m} \sigma _m-1\right){}^4},(k+1) a^2+6 \sqrt{k} \sigma _k a^2-\frac{3}{2} c^2 (m+1)-9 c^2 \sqrt{m} \sigma _m,3 (k+1) a^2+18 \sqrt{k} \sigma _k a^2-\frac{c^2 }{2}(m+1)-3 c^2 \sqrt{m} \sigma _m  $\\
 63&$ \frac{2 i a \sqrt{k} \sigma  \sigma _k \left(\sigma _m^{3/2}+\sqrt[4]{m} \text{sn}(\eta |m)\right) \left(\text{cn}(\xi |k)-\frac{\sigma _k \sqrt{-\frac{k'}{k}}}{\text{cn}(\xi |k)}\right)}{c \left(\sqrt{m}-\sigma _m\right) \sigma _m \left(\sigma _m^{3/2}-\sqrt[4]{m} \text{sn}(\eta |m)\right)} $&$  -\frac{64 \sigma _k \left(\sqrt{-k'}-\sqrt{k} \sigma _k\right){}^2 \sqrt{-k k'}}{\left(\sqrt{m} \sigma _m-1\right){}^4},(1-2 k) a^2+6 \sigma _k \sqrt{-k k'} a^2-\frac{3}{2} c^2 (m+1)-9 c^2 \sqrt{m} \sigma _m,3 (1-2 k) a^2+18\sigma _k \sqrt{-k k'} a^2-\frac{c^2}{2} (m+1)-3 c^2 \sqrt{m} \sigma _m  $\\
 64&$ \frac{2 i a \sigma  \sigma _k \left(\sigma _m^{3/2}+\sqrt[4]{m} \text{sn}(\eta |m)\right) \left(\sqrt{k'}-\text{dn}(\xi |k)^2 \sigma _k\right)}{\text{c$\sigma $}_m\left(\sqrt{m}-\sigma _m\right)\text{   }\text{dn}(\xi |k) \left(\sigma _m^{3/2}-\sqrt[4]{m} \text{sn}(\eta |m)\right)} $&$  -\frac{64 \sigma _k \left(\sqrt{k'}-\sigma _k\right){}^2 \sqrt{k'}}{\left(\sqrt{m} \sigma _m-1\right){}^4},(k-2) a^2+6 \sigma _k \sqrt{k'} a^2-\frac{3}{2} c^2 (m+1)-9 c^2 \sqrt{m} \sigma _m,3 (k-2) a^2+18 \sigma _k \sqrt{k'} a^2-\frac{1}{2} c^2 (m+1)-3 c^2 \sqrt{m} \sigma _m  $\\
 65&$ \frac{a \sigma  i\left(\sqrt{k}-\sigma _k\right) \left(\sigma _k^{3/2}+\sqrt[4]{k} \text{sn}(\xi |k)\right) \left(\sigma _m^{3/2}+\sqrt[4]{m} \text{sn}(\eta |m)\right)}{c \sigma _m \left(\sqrt{m}-\sigma _m\right) \left(\sigma _k^{3/2}-\sqrt[4]{k} \text{sn}(\xi |k)\right) \left(\sigma _m^{3/2}-\sqrt[4]{m} \text{sn}(\eta |m)\right)} $&$  \frac{\left(\sqrt{k} \sigma _k-1\right){}^4}{\left(\sqrt{m} \sigma _m-1\right){}^4},\frac{1}{2}(-k-1) a^2-3 \sqrt{k} \sigma _k a^2-\frac{3}{2} c^2 (m+1)-9c^2 \sqrt{m} \sigma _m,-\frac{3}{2} (k+1) a^2-9 \sqrt{k} \sigma _k a^2-\frac{c^2}{2} (m+1)-3 c^2 \sqrt{m} \sigma _m  $\\
 66&$ \frac{2 a \sqrt{k} \sigma  \text{sn}(\xi |k) \left(\sqrt[4]{m'} \left(-\sigma _m\right){}^{3/2}+\sqrt[4]{-m} \text{cn}(\eta |m)\right)}{c\left(\sqrt{m'} \sigma _m+\sqrt{-m}\right) \left(\sqrt[4]{-m} \text{cn}(\eta |m)-\left(-\sigma _m\right){}^{3/2} \sqrt[4]{m'}\right) } $&$  \frac{16 k}{\left(\sqrt{m'} \sigma _m+\sqrt{-m}\right){}^4},(k+1) a^2+\frac{3}{2} c^2 (2 m-1)+9 c^2 \sigma _m \sqrt{-m m'},3 (k+1) a^2+\frac{1}{2} c^2 (2 m-1)+3 c^2 \sigma _m \sqrt{-m m'}  $\\
 67&$ \frac{2 a \sigma \text{  }\sqrt{k'}\text{sc}(\xi |k) \left(\sqrt[4]{m'} \left(-\sigma _m\right){}^{3/2}+\sqrt[4]{-m} \text{cn}(\eta |m)\right)}{c\left(\sqrt{m'} \sigma _m+\sqrt{-m}\right) \left(\sqrt[4]{-m} \text{cn}(\eta |m)-\left(-\sigma _m\right){}^{3/2} \sqrt[4]{m'}\right)} $&$  \frac{16 k'}{\left(\sqrt{m'} \sigma _m+\sqrt{-m}\right){}^4},(k-2) a^2+\frac{3}{2} c^2 (2 m-1)+9 c^2 \sigma _m \sqrt{-m m'},3 (k-2) a^2+\frac{1}{2} c^2 (2 m-1)+3 c^2 \sigma _m \sqrt{-m m'}  $\\
 68&$ \frac{2 i a \sigma \text{  }\sqrt{k k'}\text{sd}(\xi |k) \left(\sqrt[4]{m'} \left(-\sigma _m\right){}^{3/2}+\sqrt[4]{-m} \text{cn}(\eta |m)\right)}{c\left(\sqrt{m'} \sigma _m+\sqrt{-m}\right) \left(\sqrt[4]{-m} \text{cn}(\eta |m)-\left(-\sigma _m\right){}^{3/2} \sqrt[4]{m'}\right)} $&$  -\frac{16 k k'}{\left(\sqrt{m'} \sigma _m+\sqrt{-m}\right){}^4},(1-2 k) a^2+\frac{3}{2} c^2 (2 m-1)+9 c^2 \sigma _m \sqrt{-m m'},(3-6 k) a^2+\frac{1}{2} c^2 (2 m-1)+3 c^2 \sigma _m \sqrt{-m m'}  $\\
 69&$ \frac{2 a k \sigma  \text{cn}(\xi |k) \text{sn}(\xi |k) }{c\left(\sqrt{m'} \sigma _m+\sqrt{-m}\right) \text{dn}(\xi |k)\text{  }}\newline\times \frac{\left(\sqrt[4]{m'} \left(-\sigma _m\right){}^{3/2}+\sqrt[4]{-m} \text{cn}(\eta |m)\right)}{\left(\sqrt[4]{-m} \text{cn}(\eta |m)-\left(-\sigma _m\right){}^{3/2} \sqrt[4]{m'}\right)} $&$  \frac{16 k^2}{\left(\sqrt{m'} \sigma _m+\sqrt{-m}\right){}^4},-2 (k-2) a^2+\frac{3}{2} c^2 (2 m-1)+9 c^2 \sigma _m \sqrt{-m m'},-6 (k-2) a^2+\frac{1}{2} c^2 (2 m-1)+3 c^2 \sigma _m \sqrt{-m m'}  $\\
 70&$ \frac{a k \sigma \sigma _k \text{cn}(\xi |k)\text{  }}{c\text{   }\left(\sqrt{m'} \sigma _m+\sqrt{-m}\right)\left(\sqrt{k'}-\text{dn}(\xi |k)\right)}\newline\times \frac{\left(\sqrt[4]{m'} \left(-\sigma _m\right){}^{3/2}+\sqrt[4]{-m} \text{cn}(\eta |m)\right)}{\left(\sqrt[4]{-m} \text{cn}(\eta |m)-\left(-\sigma _m\right){}^{3/2} \sqrt[4]{m'}\right)} $&$  \frac{k^2}{\left(\sqrt{m'} \sigma _m+\sqrt{-m}\right){}^4},-\frac{1}{2} (k-2) a^2+\frac{3}{2} c^2 (2 m-1)+9 c^2 \sigma _m \sqrt{-m m'},-\frac{3}{2} (k-2) a^2+\frac{c^2}{2} (2 m-1)+3c^2 \sigma _m \sqrt{-m m'}  $\\
 71&$ \frac{i a \sigma \sigma _k \sqrt{k'}  \text{dn}(\xi |k) }{c\text{   }\left(\sqrt{m'} \sigma _m+\sqrt{-m}\right)\left(\sqrt{k} \text{sn}(\xi |k)-1\right)}\newline\times \frac{\left(\sqrt[4]{m'} \left(-\sigma _m\right){}^{3/2}+\sqrt[4]{-m} \text{cn}(\eta |m)\right)}{\left(\sqrt[4]{-m} \text{cn}(\eta |m)-\left(-\sigma _m\right){}^{3/2} \sqrt[4]{m'}\right)} $&$  \frac{k'^2}{\left(\sqrt{m'} \sigma _m+\sqrt{-m}\right){}^4},-\frac{1}{2} (k+1) a^2+\frac{3}{2} c^2 (2 m-1)+9 c^2 \sigma _m \sqrt{-m m'},-\frac{3}{2} (k+1) a^2+\frac{c^2 }{2}(2 m-1)+3 c^2 \sigma _m \sqrt{-m m'}  $\\
 72&$ \frac{a \sigma (1-2 \mathrm{\theta} [0.5 -k]) \text{sn}(\xi |k) }{c\left(\sqrt{m'} \sigma _m+\sqrt{-m}\right) (\text{cn}(\xi |k)+1)\text{  }}\newline\times \frac{\left(\sqrt[4]{m'} \left(-\sigma _m\right){}^{3/2}+\sqrt[4]{-m} \text{cn}(\eta |m)\right)}{\left(\sqrt[4]{-m} \text{cn}(\eta |m)-\left(-\sigma _m\right){}^{3/2} \sqrt[4]{m'}\right)} $&$  \frac{1}{\left(\sqrt{m'} \sigma _m+\sqrt{-m}\right){}^4},\left(k-\frac{1}{2}\right) a^2+\frac{3}{2} c^2 (2 m-1)+9 c^2 \sigma _m \sqrt{-m m'},\left(3 k-\frac{3}{2}\right) a^2+\frac{1}{2} c^2 (2 m-1)+3 c^2 \sigma _m \sqrt{-m m'}  $\\
 73&$ \frac{2 a \sigma  \left(\sqrt{k} \sigma _k \text{sn}(\xi |k)^2+1\right) }{c \left(\sqrt{m'} \sigma _m+\sqrt{-m}\right) \text{sn}(\xi |k)}\newline\times \frac{\left(\sqrt[4]{m'} \left(-\sigma _m\right){}^{3/2}+\sqrt[4]{-m} \text{cn}(\eta |m)\right)}{\left(\sqrt[4]{-m} \text{cn}(\eta |m)-\left(-\sigma _m\right){}^{3/2} \sqrt[4]{m'}\right)} $&$  \frac{64 \sqrt{k} \sigma _k \left(\sigma _k+\sqrt{k}\right){}^2}{\left(\sqrt{m'} \sigma _m+\sqrt{-m}\right){}^4},(k+1) a^2+6 \sqrt{k} \sigma _k a^2+\frac{3}{2} c^2 (2 m-1)+9 c^2 \sigma _m \sqrt{-m m'},3 (k+1) a^2+18 \sqrt{k} \sigma _k a^2+\frac{1}{2} c^2 (2 m-1)+3 c^2 \sigma _m \sqrt{-m m'}  $\\

\hline
 74&$ \frac{2 i a \sqrt{k} \sigma  \sigma _k \left(\text{cn}(\xi |k)-\frac{\sigma _k \sqrt{-\frac{k'}{k}}}{\text{cn}(\xi |k)}\right) }{c\left(\sqrt{m'} \sigma _m+\sqrt{-m}\right) }\newline\times \frac{\left(\sqrt[4]{m'} \left(-\sigma _m\right){}^{3/2}+\sqrt[4]{-m} \text{cn}(\eta |m)\right)}{\left(\sqrt[4]{-m} \text{cn}(\eta |m)-\left(-\sigma _m\right){}^{3/2} \sqrt[4]{m'}\right)} $&$  -\frac{64 \sigma _k \left(\sqrt{-k'}-\sqrt{k} \sigma _k\right){}^2 \sqrt{-k k'}}{\left(\sqrt{m'} \sigma _m+\sqrt{-m}\right){}^4},(1-2 k) a^2+6 \sigma _k \sqrt{-k k'} a^2+\frac{3}{2} c^2 (2 m-1)+9 c^2 \sigma _m \sqrt{-m m'},(3-6 k) a^2+18 \sigma _k \sqrt{-k k'} a^2+\frac{1}{2} c^2 (2 m-1)+3 c^2 \sigma _m \sqrt{-m m'}  $\\
 75&$ \frac{2 i a \sigma  \sigma _k \left(\text{dn}(\xi |k)^2 \sigma _k-\sqrt{k'}\right) }{c \left(\sqrt{m'} \sigma _m+\sqrt{-m}\right)\text{dn}(\xi |k)}\newline\times \frac{\left(\sqrt[4]{m'} \left(-\sigma _m\right){}^{3/2}+\sqrt[4]{-m} \text{cn}(\eta |m)\right)}{\left(\sqrt[4]{-m} \text{cn}(\eta |m)-\left(-\sigma _m\right){}^{3/2} \sqrt[4]{m'}\right)} $&$  -\frac{64 \sigma _k \left(\sqrt{k'}-\sigma _k\right){}^2 \sqrt{k'}}{\left(\sqrt{m'} \sigma _m+\sqrt{-m}\right){}^4},(k-2) a^2+6 \sigma _k \sqrt{k'} a^2+\frac{3}{2} c^2 (2 m-1)+9 c^2 \sigma _m \sqrt{-m m'},3 (k-2) a^2+18 \sigma _k \sqrt{k'} a^2+\frac{1}{2} c^2 (2 m-1)+3 c^2 \sigma _m \sqrt{-m m'}  $\\
 76&$ \frac{a \sigma  \left(\sqrt{k}-\sigma _k\right) \left(\sigma _k^{3/2}+\sqrt[4]{k} \text{sn}(\xi |k)\right) }{c\left(\sqrt{m'} \sigma _m+\sqrt{-m}\right)\left(\sigma _k^{3/2}-\sqrt[4]{k} \text{sn}(\xi |k)\right)}\newline\times \frac{\left(\sqrt[4]{m'} \left(-\sigma _m\right){}^{3/2}+\sqrt[4]{-m} \text{cn}(\eta |m)\right)}{\left(\sqrt[4]{-m} \text{cn}(\eta |m)-\left(-\sigma _m\right){}^{3/2} \sqrt[4]{m'}\right)} $&$  \frac{\left(\sqrt{k} \sigma _k-1\right){}^4}{\left(\sqrt{m'} \sigma _m+\sqrt{-m}\right){}^4},\frac{1}{2}(-k-1) a^2-3 \sqrt{k} \sigma _k a^2+\frac{3}{2} c^2 (2 m-1)+9c^2 \sigma _m \sqrt{-m m'},-\frac{3}{2} (k+1) a^2-9 \sqrt{k} \sigma _k a^2+\frac{c^2}{2} (2 m-1)+3 c^2 \sigma _m \sqrt{-m m'}  $\\
 77&$ \frac{i a \sigma \text{  }\left(\sqrt{k'} \sigma _k+\sqrt{-k}\right)\left(\sqrt[4]{k'} \left(-\sigma _k\right){}^{3/2}+\sqrt[4]{-k} \text{cn}(\xi |k)\right)}{c \left(\sqrt{m'} \sigma _m+\sqrt{-m}\right) \left(\sqrt[4]{-k} \text{cn}(\xi |k)-\left(-\sigma _k\right){}^{3/2} \sqrt[4]{k'}\right)}\newline\times \frac{\left(\sqrt[4]{m'} \left(-\sigma _m\right){}^{3/2}+\sqrt[4]{-m} \text{cn}(\eta |m)\right)}{\left(\sqrt[4]{-m} \text{cn}(\eta |m)-\left(-\sigma _m\right){}^{3/2} \sqrt[4]{m'}\right)} $&$  \frac{\left(\sqrt{k'} \sigma _k+\sqrt{-k}\right){}^4}{\left(\sqrt{m'} \sigma _m+\sqrt{-m}\right){}^4},\left(k-\frac{1}{2}\right) a^2+3 \sigma _k \sqrt{-k k'} a^2+\frac{3}{2} c^2 (2 m-1)+9 c^2 \sigma _m \sqrt{-m m'},9 \sigma _k \sqrt{-k k'} a^2+\frac{1}{2}(6 k-3) a^2+\frac{c^2 }{2}(2 m-1)+3 c^2 \sigma _m \sqrt{-m m'}  $\\
 78&$ \frac{2 i a \sigma _m\sqrt{k} \sigma  \text{sn}(\xi |k)\text{  }\left(\sqrt[4]{m'} \sigma _m^{3/2}+\text{dn}(\eta |m)\right)}{c\left(\sigma _m \sqrt{m'}-1\right) \left(\text{dn}(\eta |m)-\sigma _m^{3/2} \sqrt[4]{m'}\right) } $&$  \frac{16 k}{\left(\sigma _m \sqrt{m'}-1\right){}^4},(k+1) a^2-\frac{3}{2} c^2 (m-2)+9 c^2 \sigma _m \sqrt{m'},3 (k+1) a^2-\frac{1}{2} c^2 (m-2)+3 c^2 \sigma _m \sqrt{m'}  $\\
 79&$ \frac{2 i a \sigma \text{  }\sigma _m \sqrt{k'} \text{sc}(\xi |k)\left(\sqrt[4]{m'} \sigma _m^{3/2}+\text{dn}(\eta |m)\right)}{c\left(\sigma _m \sqrt{m'}-1\right) \left(\text{dn}(\eta |m)-\sigma _m^{3/2} \sqrt[4]{m'}\right) } $&$  \frac{16 k'}{\left(\sigma _m \sqrt{m'}-1\right){}^4},(k-2) a^2-\frac{3}{2} c^2 (m-2)+9 c^2 \sigma _m \sqrt{m'},3 (k-2) a^2-\frac{1}{2} c^2 (m-2)+3 c^2 \sigma _m \sqrt{m'}  $\\
 80&$ \frac{2 a \sigma \text{  }\sigma _m \sqrt{k k'}\text{sd}(\xi |k) \left(\sqrt[4]{m'} \sigma _m^{3/2}+\text{dn}(\eta |m)\right)}{c \left(\sigma _m \sqrt{m'}-1\right)\left(\text{dn}(\eta |m)-\sigma _m^{3/2} \sqrt[4]{m'}\right) } $&$  -\frac{16 k k'}{\left(\sigma _m \sqrt{m'}-1\right){}^4},(1-2 k) a^2-\frac{3}{2} c^2 (m-2)+9 c^2 \sigma _m \sqrt{m'},(3-6 k) a^2-\frac{1}{2} c^2 (m-2)+3 c^2 \sigma _m \sqrt{m'}  $\\
 81&$ \frac{2 i a k \sigma  \sigma _m\text{  }\text{cn}(\xi |k) \text{sn}(\xi |k)\left(\sqrt[4]{m'} \sigma _m^{3/2}+\text{dn}(\eta |m)\right)}{c\left(\sigma _m \sqrt{m'}-1\right) \text{dn}(\xi |k) \left(\text{dn}(\eta |m)-\sigma _m^{3/2} \sqrt[4]{m'}\right) } $&$  \frac{16 k^2}{\left(\sigma _m \sqrt{m'}-1\right){}^4},-2 (k-2) a^2-\frac{3}{2} c^2 (m-2)+9 c^2 \sigma _m \sqrt{m'},-6 (k-2) a^2-\frac{1}{2} c^2 (m-2)+3 c^2 \sigma _m \sqrt{m'}  $\\
 82&$ \frac{i a k \sigma \text{  }\sigma _k \sigma _m \text{cn}(\xi |k)\left(\sigma _m^{3/2}\sqrt[4]{m'} +\text{dn}(\eta |m)\right)}{c \left(\sigma _m \sqrt{m'}-1\right)\left(\sqrt{k'}-\text{dn}(\xi |k)\right) \left(\text{dn}(\eta |m)-\sigma _m^{3/2} \sqrt[4]{m'}\right) } $&$ \frac{k^2}{\left(\sigma _m \sqrt{m'}-1\right){}^4},-\frac{1}{2} (k-2) a^2-\frac{3}{2} c^2 (m-2)+9 c^2 \sigma _m \sqrt{m'},-\frac{3}{2} (k-2) a^2-\frac{c^2}{2} (m-2)+3 c^2 \sigma _m \sqrt{m'}$\\
 83&$ \frac{a \sigma \text{  }\sigma _k \sigma _m \sqrt{k'}\text{dn}(\xi |k) \left(\sigma _m^{3/2}\sqrt[4]{m'} +\text{dn}(\eta |m)\right)}{c\left(\sigma _m \sqrt{m'}-1\right) \left(\sqrt{k} \text{sn}(\xi |k)-1\right) \left(\text{dn}(\eta |m)-\sigma _m^{3/2} \sqrt[4]{m'}\right) } $&$ \frac{k'^2}{\left(\sigma _m \sqrt{m'}-1\right){}^4},-\frac{1}{2} (k+1) a^2-\frac{3}{2} c^2 (m-2)+9 c^2 \sigma _m \sqrt{m'},-\frac{3}{2} (k+1) a^2-\frac{c^2}{2} (m-2)+3c^2 \sigma _m \sqrt{m'}$\\
 84&$ \frac{a \sigma \text{  }\sigma _m (1-2 \mathrm{\theta} [0.5 -k]) \text{sn}(\xi |k)}{c\text{  }\left(\sigma _m^{3/2} \sqrt[4]{m'}-1\right) \left(\sqrt[4]{m'} \sigma _m^{3/2}+1\right) (\text{cn}(\xi |k)+1)}\newline\times \frac{\left(\sigma _m^{3/2}\sqrt[4]{m'} +\text{dn}(\eta |m)\right)}{\left(\sigma _m^{3/2} \sqrt[4]{m'}-\text{dn}(\eta |m)\right)} $&$  \frac{1}{\left(\sigma _m \sqrt{m'}-1\right){}^4},\left(k-\frac{1}{2}\right) a^2-\frac{3}{2} c^2 (m-2)+9 c^2 \sigma _m \sqrt{m'},\left(3 k-\frac{3}{2}\right) a^2-\frac{1}{2} c^2 (m-2)+3 c^2 \sigma _m \sqrt{m'}  $\\
 85&$ \frac{2 i a \sigma \sigma _m \left(\sqrt{k} \sigma _k \text{sn}(\xi |k)^2+1\right)\text{  }\left(\sqrt[4]{m'} \sigma _m^{3/2}+\text{dn}(\eta |m)\right)}{c \left(\sigma _m \sqrt{m'}-1\right)\text{sn}(\xi |k) \left(\text{dn}(\eta |m)-\sigma _m^{3/2} \sqrt[4]{m'}\right) } $&$  \frac{64 \sqrt{k} \sigma _k \left(\sigma _k+\sqrt{k}\right){}^2}{\left(\sigma _m \sqrt{m'}-1\right){}^4},(k+1) a^2+6 \sqrt{k} \sigma _k a^2-\frac{3}{2} c^2 (m-2)+9 c^2 \sigma _m \sqrt{m'},3 (k+1) a^2+18 \sqrt{k} \sigma _k a^2-\frac{1}{2} c^2 (m-2)+3 c^2 \sigma _m \sqrt{m'}  $\\
 86&$ \frac{2 a\text{  }\sigma  \sigma _k \sigma _m\sqrt{k}\left(k+\sigma _k \sqrt{-k k'}\right) \left(\text{cn}(\xi |k)-\frac{\sigma _k \sqrt{-\frac{k'}{k}}}{\text{cn}(\xi |k)}\right) }{c \left(\sigma _m \sqrt{m'}-1\right)}\newline\times \frac{\left(\sigma _m^{3/2}\sqrt[4]{m'} +\text{dn}(\eta |m)\right)}{\left(\sigma _m^{3/2} \sqrt[4]{m'}-\text{dn}(\eta |m)\right)} $&$  -\frac{64 \sigma _k \left(\sqrt{-k'}-\sqrt{k} \sigma _k\right){}^2 \sqrt{-k k'}}{\left(\sigma _m \sqrt{m'}-1\right){}^4},(1-2 k) a^2+6 \sigma _k \sqrt{-k k'} a^2-\frac{3}{2} c^2 (m-2)+9 c^2 \sigma _m \sqrt{m'},(3-6 k) a^2+18 \sigma _k \sqrt{-k k'} a^2-\frac{1}{2} c^2 (m-2)+3 c^2 \sigma _m \sqrt{m'}  $\\
 87&$ \frac{2 a \sigma  \sigma _k \sigma _m \left(\sqrt{k'}-\text{dn}(\xi |k)^2 \sigma _k\right) \left(\sqrt[4]{m'} \sigma _m^{3/2}+\text{dn}(\eta |m)\right)}{c\left(\sigma _m \sqrt{m'}-1\right) \text{dn}(\xi |k) \left(\text{dn}(\eta |m)-\sigma _m^{3/2} \sqrt[4]{m'}\right) } $&$  -\frac{64 \sigma _k \left(\sqrt{k'}-\sigma _k\right){}^2 \sqrt{k'}}{\left(\sigma _m \sqrt{m'}-1\right){}^4},(k-2) a^2+6 \sigma _k \sqrt{k'} a^2-\frac{3}{2} c^2 (m-2)+9 c^2 \sigma _m \sqrt{m'},3 (k-2) a^2+18 \sigma _k \sqrt{k'} a^2-\frac{1}{2} c^2 (m-2)+3 c^2 \sigma _m \sqrt{m'}  $\\
 88&$ \frac{a \sigma  \left(\sigma _m\sqrt{k}-\sigma _k\right) \left(\sigma _k^{3/2}+\sqrt[4]{k} \text{sn}(\xi |k)\right)\text{  }}{c \left(\sigma _m \sqrt{m'}-1\right) \left(\sigma _k^{3/2}-\sqrt[4]{k} \text{sn}(\xi |k)\right) }\newline\times \frac{\left(\sigma _m^{3/2}\sqrt[4]{m'} +\text{dn}(\eta |m)\right)}{\left(\sigma _m^{3/2} \sqrt[4]{m'}-\text{dn}(\eta |m)\right)} $&$  \frac{\left(\sqrt{k} \sigma _k-1\right){}^4}{\left(\sigma _m \sqrt{m'}-1\right){}^4},\frac{1}{2}(-k-1) a^2-3 \sqrt{k} \sigma _k a^2-\frac{3}{2} c^2 (m-2)+9 c^2 \sigma _m \sqrt{m'},-\frac{3}{2} (k+1) a^2-9 \sqrt{k} \sigma _k a^2-\frac{1}{2} c^2 (m-2)+3 c^2 \sigma _m \sqrt{m'}  $\\

\hline
 89&$ \frac{a \sigma  \sigma _m\text{  }\left(\sqrt{k'} \sigma _k+\sqrt{-k}\right) \left(\sqrt[4]{k'} \left(-\sigma _k\right){}^{3/2}+\sqrt[4]{-k} \text{cn}(\xi |k)\right)}{c \left(\sigma _m \sqrt{m'}-1\right)\left(\sqrt[4]{-k} \text{cn}(\xi |k)-\left(-\sigma _k\right){}^{3/2} \sqrt[4]{k'}\right)}\newline\times \frac{\left(\sigma _m^{3/2}\sqrt[4]{m'} +\text{dn}(\eta |m)\right)}{\left(\sigma _m^{3/2} \sqrt[4]{m'}-\text{dn}(\eta |m)\right)} $&$  \frac{\left(\sqrt{k'} \sigma _k+\sqrt{-k}\right){}^4}{\left(\sigma _m \sqrt{m'}-1\right){}^4},\left(k-\frac{1}{2}\right) a^2+3 \sigma _k \sqrt{-k k'} a^2-\frac{3}{2} c^2 (m-2)+9 c^2 \sigma _m \sqrt{m'},\left(3 k-\frac{3}{2}\right) a^2+9 \sigma _k \sqrt{-k k'} a^2-\frac{1}{2} c^2 (m-2)+3 c^2 \sigma _m \sqrt{m'}  $\\
 90&$ \frac{i a \sigma  \sigma _m \left(\sqrt{k'} \sigma _k+\sqrt{-k}\right) \left(\sqrt[4]{k'} \sigma _k^{3/2}+\text{dn}(\xi |k)\right)}{c \left(\sigma _m \sqrt{m'}-1\right) \left(\text{dn}(\xi |k)-\sigma _k^{3/2} \sqrt[4]{k'}\right)}\newline\times \frac{\left(\sigma _m^{3/2}\sqrt[4]{m'} +\text{dn}(\eta |m)\right)}{\left(\sigma _m^{3/2} \sqrt[4]{m'}-\text{dn}(\eta |m)\right)} $&$  \frac{\left(\sigma _k \sqrt{k'}-1\right){}^4}{\left(\sigma _m \sqrt{m'}-1\right){}^4},-\frac{1}{2} (k-2) a^2+3 \sigma _k \sqrt{k'} a^2-\frac{3}{2} c^2 (m-2)+9 c^2 \sigma _m \sqrt{m'},-\frac{3}{2} (k-2) a^2+9 \sigma _k \sqrt{k'} a^2-\frac{1}{2} c^2 (m-2)+3 c^2 \sigma _m \sqrt{m'} $\\
\end{longtable} \end{center}
\small Note: $  \mathrm{\theta}[x]$ is $\mathbf{\theta}$ function,it satisfies $x< 0,~\mathrm{\theta} [x]=0; x\geq 0,~\mathrm{\theta} [x]=1.~  \sigma _m^{2}=1,~\sigma _k^{2}=1$.
\section{Conclusion}

So we can conclude that If two Jacobi elliptic function solutions of a NDE satisfy  (1) and (2),~the two solutions are equivalent. According to the usage of the modulus and phase transformation of elliptic functions to mKdV,~a lot of the eqivalent Jacobi elliptic function solutions can be obtained when we knowing one of the Jacobi elliptic function solutions to some NDEs.~It can be used  to test the eqivalent of the other  NDEs' Jacobi elliptic function solutions.

\section*{Appendix: Definition and Jacobi elliptic function}

Common  Jacobi elliptic function have three kinds$^{[6,7]}$
:$ \text{sn}(x\vert k)$,~$\text{cn}(x|
k)$ and $\text{dn}(x| k)$,
where $k$ is the module of the Jacobi elliptic function.
They meet the relations
\begin{eqnarray}
\label{eq1}
\text{sn}(x| k)^2+\text{cn}(x| k)^2= 1,
\text{dn}(x| k)^2+k
\text{sn}(x| k)^2= 1,~\text{dn}(x| k)^2-k \text{cn}(x| k)^2=
1-k={k}'.
\end{eqnarray}
It is known that
\begin{eqnarray}
\label{eq2}
\text{sn}(-x| k) =-\text{sn}(x| k),~\text{cn}(-x| k) =
\text{cn}(x| k),~\text{dn}(-x| k) = \text{dn}(x| k),
\end{eqnarray}
where ${k}'$ is complementary  module of Jacobi elliptic function.
Except $\text{sn}(x\vert k)$,~$\text{cn}(x|
k)$ and $\text{dn}(x| k)$,~other nine Jacobi elliptic functions can be defined
\\
$\begin{array}{ccccccc}
 1/\text{sn}(x| k)\equiv \text{ns}(x| k),~\text{dn}(x|k)/\text{cn}(x| k)\equiv \text{dc}(x| k),~\text{cn}(x|k)/\text{dn}(x| k)\equiv \text{cd}(x| k),~\\
 1/\text{cn}(x| k)\equiv \text{nc}(x| k),~\text{dn}(x|k)/\text{sn}(x| k)\equiv \text{ds}(x| k),~\text{sn}(x|k)/\text{dn}(x| k)\equiv \text{sd}(x| k),~\\
 1/\text{dn}(x| k)\equiv \text{nd}(x| k),~\text{sn}(x|k)/\text{cn}(x| k)\equiv \text{sc}(x| k),~\text{cn}(x|k)/\text{sn}(x| k))\equiv \text{cs}(x| k).~\\
 \end{array}$\\\\
Jacobi elliptic function has some common phase transformation,~such as
\begin{eqnarray}
\text{sn}(x+2pK+2iq{K}'|k)=
 (-1)^{p}\text{sn}(x|k),
 \end{eqnarray}
\begin{eqnarray}
\text{cn}(x+2pK+2iq{K}'| k)= (-1)^{p+q}\text{cn}(x| k),~\end{eqnarray}
\begin{eqnarray}
\text{dn}(x+2pK+2iq{K}'| k)= (-1)^q\text{dn}(x| k),~\end{eqnarray}
\begin{eqnarray}
\text{sn}(x+K+i{K}'| k)=\frac{\text{dc}(x| k)}{\sqrt k },~\end{eqnarray}
\begin{eqnarray}
\text{cn}(x+K+i{K}'| k)=-i\sqrt {{{k}'} \mathord{\left/ {\vphantom
{{{k}'} k}} \right.~\kern-\nulldelimiterspace} k} \text{nc}(x| k),~\end{eqnarray}
\begin{eqnarray} \text{dn}(x+K+i{K}'|
k)=-i\sqrt {{k}'} \text{sc}(x| k),~\quad _{ } \end{eqnarray}
\begin{eqnarray}
\text{sn}(x+K| k)=\text{cd}(x| k),~\end{eqnarray}
\begin{eqnarray}
\text{cn}(x+K| k)=\sqrt {{k}'} \text{sd}(x| k),~\end{eqnarray}
\begin{eqnarray}
\text{dn}(x+K| k)=\sqrt {{k}'} \text{nd}(x| k),\end{eqnarray}
\begin{eqnarray}
\text{sn}(x+i{K}'| k)=\frac{1}{\sqrt k }\text{ns}(x| k),\end{eqnarray}
\begin{eqnarray}
\text{cn}(x+i{K}'| k)=\frac{-i}{\sqrt k }\text{ds}(x| k),\end{eqnarray}
\begin{eqnarray}
\text{dn}(x+i{K}'| k)=\text{cs}(x| k),
\end{eqnarray}
where $i^2=-1,~K=\int_0^{\pi /2} {[1-k\sin ^2(\theta )]} ^{-\frac{1}{2}}d\theta
,~{K}'=\int_0^{\pi /2} {[1-{k}'\sin ^2(\theta )]} ^{-\frac{1}{2}}d\theta $,
$p$ and $q$ are integers.

\end{document}